\documentclass{amsart}

\usepackage{amsfonts,amssymb,mathrsfs}
\usepackage{graphics,color}
\newtheorem{definition}{Definition}
\newtheorem{corollary}{Corollary}
 
\newtheorem{lemma}{Lemma} 
\newtheorem{theorem}{Theorem} 
 
\theoremstyle{remark}
\newtheorem{example}{Example}
\newtheorem{remark}{Remark} 
\title[Max-plus $(A,B)$-invariant spaces and control of timed discrete event systems]{Max-plus $(A,B)$-invariant spaces and control of timed discrete event systems}
\author{Ricardo David Katz}
\address{CONICET. Postal address: Instituto de Matem\'atica ``Beppo Levi'', 
Universidad Nacional de Rosario, Av. Pellegrini 250, 2000 Rosario, Argentina.}
\email{rkatz@fceia.unr.edu.ar}
\thanks{This work was partially supported by INRIA}

\DeclareMathAlphabet{\mathbbold}{U}{bbold}{m}{n}

\newcommand{\set}[2]{\{#1\mid\,#2\}}
\newcommand{\R}{\mathbb{R}}
\newcommand{\N}{\mathbb{N}}
\newcommand{\Z}{\mathbb{Z}}

\newcommand{\Zm}{\Z\cup\{-\infty\}}
\newcommand{\Rm}{\R\cup\{-\infty\}}

\newcommand{\Np}{\N\cup\{+\infty\}}

\newcommand{\rmax}{\R_{\max}}
\newcommand{\rmaxb}{\overline{\R}_{\max}}

\newcommand{\rminb}{\overline{\R}_{\min}}

\newcommand{\nmin}{\N_{\min}}

\newcommand{\zmin}{\Z_{\min}}
\newcommand{\zminb}{\overline{\Z}_{\min}}
\newcommand{\zmax}{\Z_{\max}}
\newcommand{\zmaxb}{\overline{\Z}_{\max}}

\newcommand{\cG}{\mathcal{G}}
\newcommand{\cB}{\mathcal{B}}

\newcommand{\cS}{\mathcal{S}}

\newcommand{\cX}{\mathcal{X}}
\newcommand{\cY}{\mathcal{Y}}
\newcommand{\cZ}{\mathcal{Z}}
\newcommand{\cK}{\mathcal{K}}
\newcommand{\im}{\mbox{\rm Im}\,}
\newcommand{\vol}{\mbox{\rm vol}\,}
\newcommand{\card}{\mbox{\rm card}\,}
\newcommand{\smallvect}{\mbox{\rm span}\,}

\begin{document}

\begin{abstract}
The concept of $(A,B)$-invariant subspace (or controlled invariant) of 
a linear dynamical system is extended to linear systems over the max-plus 
semiring. Although this extension presents several difficulties, which are 
similar to those encountered in the same kind of extension to linear 
dynamical systems over rings, it appears capable of providing solutions to 
many control problems like in the cases of linear systems over fields or rings. 
Sufficient conditions are given for computing the maximal $(A,B)$-invariant 
subspace contained in a given space and the existence of linear state 
feedbacks is discussed. An application to the study of transportation 
networks which evolve according to a timetable is considered.
\end{abstract}

\maketitle

\section{Introduction}

The geometric approach to the theory of linear 
dynamical systems has provided deep insights 
and elegant solutions to many control problems, 
such as the disturbance decoupling problem, the block 
decoupling problem, and the model matching problem 
(see~\cite{wonham} and the references therein). 
The concept of $(A,B)$-invariant subspace 
(or controlled invariant subspace, see~\cite{BasMar91}) 
has played a significant role in the development of this approach. 

It is natural to try to apply the same kind of methods to discrete 
event systems. Several mathematical models have been proposed, see 
in particular~\cite{CassLafoOlsd} for a survey of the following 
approaches. Ramadge and Wonham~\cite{ramadge87a} initiated the logical, 
language-theoretic approach, in which the precise ordering of the events 
is of interest and time does not play an explicit role. 
This theory addresses the synthesis of controllers in 
order to satisfy some qualitative specifications on the admissible orderings 
of the events. Another approach is the max-plus algebra based control approach 
initiated by Cohen et al.~\cite{cohen85a}, in which in addition to 
the ordering, the timing of the events plays an essential role. A third 
approach is the perturbation analysis of Cassandras 
and Ho~\cite{CassaHo83}, which deals with stochastic timed 
discrete event systems. 

The max-plus semiring is the set $\Rm$, equipped with 
$\max$ as addition and the usual sum as multiplication. 
Linear dynamical systems with coefficients in the max-plus 
semiring turn out to be useful for modeling and analyzing 
many discrete event dynamic systems subject to synchronization 
constraints (see~\cite{bcoq}). Among these, we can mention 
some manufacturing systems (Cohen et al.~\cite{cohen85a}), 
computer networks (Le Boudec and Thiran~\cite{leboudec}) and 
transportation networks (Olsder et al.~\cite{OlsSubGett98}, 
Braker~\cite{braker91,braker}, and de Vries et al.~\cite{deVDeSdeM98}). 
Many results from linear system theory have been extended to 
systems with coefficients in the max-plus semiring, such as 
the connection between spectral theory and stability questions 
(see~\cite{cohen89a}) or transfer series methods (see~\cite{bcoq}). 
Several interesting control problems have also been studied by, 
for example, Boimond et al.~\cite{BoiCott99,BoiHar00}, 
Cottenceau et al.~\cite{CottHar03} and Lhommeau~\cite{Lhommeau}. 
In contrast to the approach presented here, which is based on state space 
representation, their approach uses transfer series and residuation methods 
and therefore deals with different types of specifications. 

This motivates the attempt to extend the geometric approach, and in 
particular the concept of $(A,B)$-invariant subspace, to the 
theory of linear dynamical systems over the max-plus 
semiring, a question which is raised in~\cite{ccggq99}. 
The same kind of generalization, which was initiated by Hautus, Conte and 
Perdon, has been widely studied for linear dynamical systems over rings 
(see~\cite{hautus82,hautus84,conte94,conte95,assan,AssLafPer}).
In this paper we will see that the extension of the 
geometric approach to linear systems over the max-plus 
semiring presents similar difficulties to those 
encountered in dealing with coefficients in a ring 
rather than coefficients in a field. The $(A,B)$-invariance problem has been studied 
in the framework of formal series over some complete idempotent 
semirings by Klimann~\cite{klimann99}. 

\begin{figure}
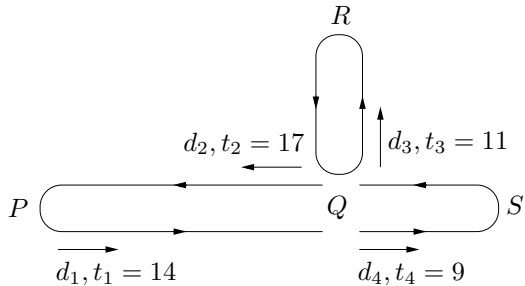

\begin{center}
\input figure1V2
\end{center}
\caption{A simple transportation network}
\label{figure1}
\end{figure}

To illustrate one of the possible applications of the results presented in 
this paper, we apply the methods presented here to the study of transportation 
networks which evolve according to a timetable. Max-plus linear models for 
transportation networks have been studied by several authors, see for 
example~\cite{OlsSubGett98,braker91,braker,deVDeSdeM98}. Let us 
consider the simple railway network given in Figure~\ref{figure1}, which 
has been borrowed from~\cite{deVDeSdeM98}. 
In this network, we assume that in the initial state there is a train 
running along each of the tracks which connect the following stations: 
$P$ with $Q$, $Q$ with $P$, $Q$ with $Q$ via $R$ and finally $Q$ with 
$Q$ via $S$. In Figure~\ref{figure1}, these tracks are denoted by 
$d_1$, $d_2$, $d_3$ and $d_4$ respectively. The traveling time on 
track $d_i$ is given by $t_i$, for $i=1,\ldots ,4$. We will assume that the following 
conditions are satisfied. A first condition is that at station $Q$ the 
trains coming from stations $P$ and $S$ have to ensure a connection to the 
train which leaves for destination $R$ and vice versa. The second 
condition is that a train cannot leave before its scheduled departure 
time which is given by a timetable. If we assume that a train leaves 
as soon as all the previous conditions have been satisfied, then the 
evolution of the transportation network can be described by a max-plus 
linear dynamical system where the scheduled departure times can be 
seen as controls (see Section~\ref{aplicacionSec}). We will see that 
the tools presented in this paper can be used to analyze this kind of network. 
For example, it is possible to determine whether there exists a 
timetable that satisfies such conditions as the following. 
A first condition could be that the time between two consecutive 
departures of trains in the same direction be less than a certain 
given bound. As a second condition we could require that the time 
that passengers have to wait to make some connections be less 
than another given bound. Of course, more general specifications could be 
analyzed. We show how to compute a 
timetable which satisfies these requirements when it exists. 
For instance, suppose that in the railway network given in 
Figure~\ref{figure1} we want the time between two 
consecutive departures of trains in the same direction to be less 
than $15$ time units and the maximal time that passengers 
have to wait to make any connection to be less than $4$ time units. 
In Section~\ref{aplicacionSec} we show that this is possible 
and give a timetable which satisfies these requirements.

This paper is organized as follows. In Section~\ref{geomABinvSec}, 
after a short introduction to max-plus type semirings, we introduce 
the concept of geometrically $(A,B)$-invariant semimodule and generalize 
the Wonham fixed point algorithm (which is used to compute the maximal 
$(A,B)$-invariant subspace contained in a given space, see~\cite{wonham}) 
to max-plus algebra. In Section~\ref{volumeSec} we introduce the concept of volume of 
a semimodule and study its properties. In Section~\ref{finitevolumeSec} 
we use volume arguments to show that the fixed point algorithm introduced 
in Section~\ref{geomABinvSec} converges in a finite number of steps for an 
important class of semimodules. In Section~\ref{algABinvSec} we consider the 
concept of algebraically $(A,B)$-invariant semimodule and give a method to 
decide whether a finitely generated semimodule is algebraically 
$(A,B)$-invariant. Finally, in Section~\ref{aplicacionSec} we apply the 
methods given in this paper to the study of transportation networks 
which evolve according to a timetable.

Let us finally mention that some of the results presented here were 
announced in~\cite{gk03} and considered in~\cite{katz}. 

\medskip\noindent{\em Acknowledgment.}\/
The author would like to thank S. Gaubert for many helpful 
suggestions and comments on preliminary versions of this manuscript 
and J.-J. Loiseau for useful references. He would also like to thank 
J. E. Cury and the anonymous reviewers who helped to improve 
this paper.    

\section{Geometrically $(A,B)$-invariant semimodules}\label{geomABinvSec}

Let us first recall some definitions and results. A {\em monoid} 
is a set equipped with an associative internal composition law 
which has a (two sided) neutral element. A {\em semiring} is a 
set ${\cS} $ equipped with two internal composition laws $\oplus $ 
and $\otimes $, called addition and multiplication respectively, 
such that ${\cS}$ is a commutative monoid for addition, $\cS$ is a 
monoid for multiplication, multiplication distributes over addition, 
and the neutral element for addition is absorbing for multiplication. 
We will sometimes denote by $({\cS},\oplus,\otimes,\varepsilon,e)$ 
the semiring ${\cS}$, where $\varepsilon$ and $e$ represent the 
neutral elements for addition and for multiplication respectively. We 
say that a semiring ${\cS}$ is {\em idempotent} if $x\oplus x=x$ 
for all $x\in \cS$. In this paper, we are mostly interested in 
some variants of the max-plus semiring $\rmax $, which is the set 
$\Rm$ equipped with $\oplus =\max $ and $\otimes =+$ 
(see~\cite{pin95} for an overview). Some of these variants can be 
obtained by noting that a semiring $M_{\max }$, whose set of elements 
is $M\cup\{-\infty\}$ and laws are $\oplus =\max $ and $\otimes =+ $, is 
associated with a submonoid $(M,+)$ of $(\R,+)$. Symmetrically, 
we can consider the semiring $M_{\min }$ with the set 
of elements $M\cup\{+\infty\}$ and laws $\oplus =\min $ 
and $\otimes =+$. 
For instance, taking $M=\Z$ we get the semiring $\zmax =(\Zm,\max,+)$, 
which is the main semiring we are going to work 
with, and taking $M=\N$ we get the semiring $\nmin =(\Np,\min,+)$, 
which is known as the {\em tropical semiring} (see~\cite{pin95}). 
Recall that an idempotent semiring 
$({\cS},\oplus,\otimes)$ is equipped with the {\em natural order}: 
$x\preceq y \iff x\oplus y=y$ (see for example~\cite{bcoq}). 
Sometimes it is useful to add a maximal element for the natural order  
to the semirings $M_{\max }$ and $M_{\min }$, obtaining in 
this way the {\em complete} semirings 
${\overline{M}}_{\max}=(M\cup\{\pm\infty\},\max,+)$ and 
${\overline{M}}_{\min}=(M\cup\{\pm\infty\},\min,+)$, respectively. 
Note that, in the semirings ${\overline{M}}_{\max}$ 
and ${\overline{M}}_{\min}$, the value of 
$(-\infty)+(+\infty)=(+\infty)+(-\infty)$ is determined 
by the fact that the neutral element for addition 
is absorbing for multiplication. Then, we know that 
$(-\infty)+(+\infty)=(+\infty)+(-\infty)=-\infty$ in 
${\overline{M}}_{\max}$ and 
$(-\infty)+(+\infty)=(+\infty)+(-\infty)=+\infty$ in 
${\overline{M}}_{\min}$.

We next introduce the concept of semimodules which is the analogous 
over semirings of vector spaces (we refer the reader to~\cite{GargKumar95} 
and~\cite{gaubert98n} for more details on semimodules). A (left) {\em semimodule} 
over a semiring $({\cS},\oplus,\otimes,\varepsilon_{{\cS}},e)$ 
is a commutative monoid $({\cX},\hat{\oplus })$, 
with neutral element $\varepsilon_{{\cX}}$, equipped with a map 
${\cS}\times {\cX}\to {\cX}$, $(\lambda,x) \to \lambda \cdot x$ 
(left action), which satisfies: 
\begin{eqnarray*}
(\lambda \otimes \mu)\cdot x= \lambda \cdot(\mu \cdot x)\; , \\
\lambda \cdot (x\; \hat{\oplus }\;y) =\lambda \cdot x \; \hat{\oplus }\; \lambda \cdot y\;,\\
(\lambda \oplus \mu)\cdot x = \lambda \cdot x \; \hat{ \oplus }\; \mu \cdot x \; , \\
\varepsilon_{\cS}\cdot x = \varepsilon_{\cX}\; , \\
\lambda \cdot \varepsilon_{\cX} = \varepsilon_{\cX} \; , \\
e \cdot x =x \; , 
\end{eqnarray*}
for all $x,y\in {\cX}$ and $\lambda,\mu\in {\cS}$. 
We will usually use concatenation to denote both the 
multiplication of ${\cS}$ and the left action, 
and we will denote by $\varepsilon$ both the zero 
element $\varepsilon_{{\cS}}$ of ${\cS}$ and the zero 
element $\varepsilon_{{\cX}}$ of ${\cX}$. A {\em subsemimodule} 
of ${\cX}$ is a subset ${\cZ}\subset {\cX}$ such that 
$\lambda x  \hat{\oplus } \mu y \in {\cZ}$, 
for all $x,y\in {\cZ}$ and $\lambda,\mu\in {\cS}$. 
In this paper, we will mostly consider subsemimodules 
of the {\em free semimodule} ${\cS}^n$, which is the set of 
$n$-dimensional vectors over ${\cS}$, equipped with 
the internal law $(x\hat{\oplus }y)_i=x_i\oplus y_i$ and the 
left action $(\lambda\cdot  x)_i=\lambda \otimes x_i$. 
If $G\subset {\cX}$, we will denote by $\smallvect G $ 
the subsemimodule of ${\cX}$ generated by $G$, 
that is, the set of all $x\in {\cX}$ for which there 
exists a finite number of elements $u_1,\ldots ,u_k$ 
of $G$ and a finite number of scalars 
$\lambda_1,\ldots ,\lambda_k\in {\cS}$, such that 
$x=\hat{\bigoplus}_{i=1,\ldots , k}\lambda_i u_i$. 
Finally, if $C\in {\cS}^{n\times r}$, we will denote by $\im C$ the 
subsemimodule of ${\cS}^n$ generated by the columns of $C$.

Let $({\cS},\oplus ,\otimes)$ denote a semiring. By a {\em system with 
coefficients in ${\cS}$}, or a {\em system over ${\cS}$}, we mean a 
linear dynamical system whose evolution is 
determined by a set of equations of the form
\begin{equation}\label{dynamicsystem}
x(k)=Ax(k-1)\oplus Bu(k)\; ,
\end{equation}
where $A\in \cS^{n\times n}$, $B\in \cS^{n\times q}$, and 
$x(k)\in \cS^{n\times 1}$, $u(k)\in \cS^{q\times 1}$, $k=1,2,\ldots $ 
are the sequences of state and control vectors respectively. 

We are interested in studying the following problem: 
Given a certain specification for the state space of 
system~\eqref{dynamicsystem}, which we suppose is given by a 
semimodule $\cK \subset \cS^n$, we want to compute the maximal 
set of initial states $\cK^*$ for which there exists a sequence 
of control vectors which makes the state of system~\eqref{dynamicsystem} 
stay in $\cK$ forever, that is, such that $x(k)\in \cK$ for all $k\geq 0$. 
To treat this problem it is convenient to make the following definition. 

\begin{definition}\label{ABinvariante2}
Given the matrices $A\in{\cS}^{n\times n}$ and $B\in {\cS}^{n\times q}$, 
we say that a  semimodule ${\cX} \subset {\cS}^n$ is 
{\rm (geometrically) $(A,B)$-invariant} if for all $x\in \cX$ there exists 
$u\in {\cS}^q$ such that $Ax \oplus Bu$ belongs to $\cX$.
\end{definition}
 
The proof of the following lemma is identical to the case of linear dynamical systems over rings. 
We include it for completeness. 

\begin{lemma}\label{obs1}
If ${\cK} \subset {\cS}^n$ is a semimodule, then ${\cK}^*$ is the maximal 
(geometrically) $(A,B)$-invariant semimodule contained in ${\cK}$. 
\end{lemma}

\begin{proof}
In the first place, note that a semimodule $\cX \subset {\cS}^n$ is (geometrically) 
$(A,B)$-invariant if and only if for each $x\in {\cX}$ there exists a sequence of control 
vectors such that the trajectory of the dynamical system~\eqref{dynamicsystem}, 
associated with this control sequence and the initial condition $x(0)=x$, 
is completely contained in ${\cX}$. Therefore, any (geometrically) $(A,B)$-invariant 
semimodule contained in ${\cK}$ is also contained in $\cK^*$. In the second place, 
note that $\cK^*$ is a subsemimodule of $\cS^n$ since system~\eqref{dynamicsystem} 
is linear and $\cK$ is a semimodule. Then, to prove the lemma, 
it only remains to show that $\cK^*$ is (geometrically) $(A,B)$-invariant. 
Let $x$ be an arbitrary element of $\cK^*$. We must see that there is a control
$u(1)\in \cS^q$ such that $x(1)=Ax \oplus Bu(1)$ belongs to $\cK^*$.
Since $x\in \cK^*$, we know that there exists a sequence of control
vectors $u(k)$, $k=1, 2, \ldots $, such that the trajectory
$x(0)$, $x(1)$, $x(2)$, $\ldots$ of system~\eqref{dynamicsystem},
associated with this control sequence and the initial condition $x(0)=x$,
is completely contained in $\cK$. Therefore, $x(1)\in \cK^*$
since there exists a sequence of control vectors
($u'(k)=u(k+1)$, $k=1, 2, \ldots $) which makes the state of
system~\eqref{dynamicsystem} stay in $\cK$ forever when the initial
state is $x(1)$.
\end{proof}

To tackle the previous problem in the case of max-plus type semirings, 
we generalize the classical fixed point algorithm which is 
used to compute the maximal $(A,B)$-invariant subspace contained 
in a given space (see~\cite{wonham}). With this purpose in mind, 
we set ${\cB}=\im B$ and consider the self-map $\varphi$ 
of the set of subsemimodules of ${\cS}^n$, given by: 
\begin{equation}\label{definicionphi}
\varphi({\cX})={\cX} \cap A^{-1}({\cX} \ominus {\cB})
\enspace ,
\end{equation}
where $A^{-1}({\cY})=\set{u\in {\cS}^n}{Au\in{\cY}}$ and  
${\cZ}\ominus {\cY}=\set{u\in {\cS}^n}{\exists y\in{\cY}, u\oplus y\in {\cZ}}$ 
for all ${\cZ},{\cY}\subset{\cS}^n$.
 
\begin{remark}\label{ObsComputo}
Note that when ${\cS}=\zmax$ or ${\cS}=\nmin$, if the semimodule ${\cX}$ is finitely 
generated, then the semimodule $\varphi({\cX})$ is also finitely generated. In fact, 
given the sets of generators of some finitely generated semimodules ${\cZ}$ and ${\cY}$, 
the semimodules ${\cY} \ominus {\cZ}$, $A^{-1}({\cY})$ and ${\cY} \cap {\cZ}$ 
can be expressed as the images by suitable matrices of the sets of solutions  
of appropriate max-plus linear systems of the form $Dx=Cx$ 
(see~\cite{gaubert98n} for details). Therefore, their sets of generators 
can be explicitly computed using a general elimination algorithm 
due to Butkovi\v{c} and Heged\"{u}s~\cite{butkovicH} and Gaubert~\cite{gaubert92a}. 
Then, when ${\cX}$ is finitely generated, the set of generators of 
$\varphi({\cX})$ can also be computed using this algorithm. More generally, if ${\cX}$ 
belongs to the class of rational semimodules (this class, which extends the 
notion of finitely generated semimodule, turns out to be useful in the geometric 
approach to discrete event systems, see~\cite{gk02a}), 
then $\varphi({\cX})$ is also a rational semimodule and can be computed by 
Theorem~3.5 of~\cite{gk02a}. 
\end{remark}

\begin{lemma}\label{obs2}\
A semimodule ${\cX} \subset {\cS}^n$ is (geometrically) $(A,B)$-invariant 
if and only if ${\cX}=\varphi({\cX})$.
\end{lemma}
\begin{proof}
Since
\begin{eqnarray*}
A^{-1}({\cX} \ominus {\cB}) & = &\set{x\in {\cS}^n}{Ax\in{\cX} \ominus {\cB}}= \\
& = & \set{x\in {\cS}^n}{\exists b\in {\cB}, Ax\oplus b\in {\cX}}= \\
& = & \set{x\in {\cS}^n}{\exists u\in {\cS}^q, Ax\oplus Bu\in {\cX}} \; ,
\end{eqnarray*}
we see that $A^{-1}({\cX} \ominus {\cB})$ is the set of initial 
states $x(0)$ of the dynamical system~\eqref{dynamicsystem} 
for which there exists a control $u(1)$ which makes the new 
state of the system, that is $x(1)=Ax(0)\oplus Bu(1)$, belong 
to ${\cX}$. Then, it readily follows from Definition~\ref{ABinvariante2} 
that a semimodule ${\cX} \subset {\cS}^n$ is (geometrically) $(A,B)$-invariant 
if and only if ${\cX}\subset A^{-1}({\cX} \ominus {\cB})$. Therefore, a semimodule 
${\cX} \subset {\cS}^n$ is (geometrically) $(A,B)$-invariant if and only if 
${\cX}=\varphi({\cX})$, that is, (geometrically) $(A,B)$-invariant semimodules 
are precisely the fixed points of the map $\varphi$ defined 
by~\eqref{definicionphi}.
\end{proof}

Inspired by the algorithm in the classical case, 
we define the following sequence of semimodules: 
\begin{equation}\label{algoABinv}
{\cX}_1={\cK}\; , \quad {\cX}_{r+1}=\varphi({\cX}_r)\;, \quad \forall r\in \N.
\end{equation}
Then we have the following lemma.
\begin{lemma}\label{lemaalgoAB}
Let ${\cK} \subset {\cS}^n$ be an arbitrary semimodule. 
Then the sequence of semimodules $\{{\cX}_r\}_{r\in \N}$ 
defined by~\eqref{algoABinv} is decreasing, i.e. 
${\cX}_{r+1}\subset {\cX}_r$ for all $r\in \N$. Moreover, 
if we define ${\cX}_{\omega}=\cap_{r\in \N} {\cX}_r$, then every 
(geometrically) $(A,B)$-invariant semimodule 
contained in ${\cK}$ is also contained in ${\cX}_{\omega}$. 
In particular, it follows that ${\cK}^* \subset {\cX}_{\omega}$.
\end{lemma}
\begin{proof}
The fact that the sequence of semimodules $\{{\cX}_r\}_{r\in \N}$ 
is decreasing is a consequence of the definition of the map $\varphi$: 
\[
{\cX}_{r+1}=\varphi({\cX}_r)={\cX}_r \cap A^{-1}({\cX}_r \ominus {\cB})\subset {\cX}_r, 
\]
for all $r\in \N$. 

To prove the second part of Lemma~\ref{lemaalgoAB}, 
firstly it is convenient to notice that $\varphi$ 
satisfies the following property: 
\[
\forall {\cZ},{\cY}\subset{\cS}^n\;,\enspace  
{\cZ}\subset {\cY} \Rightarrow \varphi({\cZ})\subset \varphi({\cY})\;,
\] 
that is, $\varphi $ is monotonic when the set of 
subsemimodules of ${\cS}^n$ is equipped with the order: 
${\cZ} \leq {\cY}$ if and only if ${\cZ}\subset {\cY}$. 

Now let ${\cX} \subset {\cK}$ be an arbitrary (geometrically) 
$(A,B)$-invariant semimodule. We will prove by induction 
on $r$ that ${\cX}\subset {\cX}_r$ for all $r\in \N$, 
and therefore that ${\cX}\subset \cap_{r\in \N} {\cX}_r={\cX}_{\omega}$. 
In the first place, we know that ${\cX} \subset {\cK} ={\cX}_1$. 
Since ${\cX}$ is a (geometrically) $(A,B)$-invariant semimodule, 
thanks to Lemma~\ref{obs2}, it follows that ${\cX} =\varphi({\cX})$. 
If we now assume that ${\cX} \subset {\cX}_t$, then we have:
\[
{\cX} = \varphi({\cX})\subset \varphi({\cX}_t)={\cX}_{t+1}\;.
\]
Therefore, ${\cX}\subset {\cX}_r$ for all $r\in \N$, as we wanted to show.
\end{proof}

Note that if the sequence $\{{\cX}_r\}_{r\in \N}$ 
stabilizes\footnote{Throughout this paper, we will use the word 
``stabilize'' to mean ``converge in a finite number of steps''.}, 
that is, if there exists $k\in \N$ such that ${\cX}_{k+1}={\cX}_k$, 
then our problem will be solved. Indeed, if there exists $k\in \N$ 
such that ${\cX}_k={\cX}_{k+1}=\varphi({\cX}_k)$ then, 
thanks to Lemma~\ref{obs2}, we know that ${\cX}_k$ is a 
(geometrically) $(A,B)$-invariant semimodule which is contained in ${\cK}$ 
(since ${\cX}_1={\cK}$ and by Lemma~\ref{lemaalgoAB} 
the sequence $\{{\cX}_r\}_{r\in \N}$ is decreasing). Therefore 
${\cX}_k\subset {\cK}^*$, and as by Lemma~\ref{lemaalgoAB} we know 
that ${\cK}^*\subset {\cX}_k$, it follows finally that ${\cK}^*={\cX}_k$.

\begin{example}\label{ejemplo1}
Let ${\cS}=\zmax$. Let us consider the matrices
\[
A=
\begin{pmatrix} -\infty & 0 \\ 0 & -\infty \end{pmatrix}
\enspace \mbox{ and } \enspace
B=\begin{pmatrix} 0 \\ 0\end{pmatrix}\; ,
\]
and the semimodule ${\cK}=\set{(x,y)^T\in \zmax^2}{y\geq x+1}$. 
Let us compute, in this particular case, 
the sequence of semimodules $\{{\cX}_r\}_{r\in \N}$ defined 
by~\eqref{algoABinv}. By definition we know that 
${\cX}_1={\cK}=\set{(x,y)^T\in \zmax^2}{y\geq x+1}$. Since 
there exists $\lambda \in \zmax$ such that 
$\max(y,\lambda)\geq\max(x,\lambda)+1$ (that is, 
there exists $(\lambda,\lambda)^T\in \cB$ such that 
$(x,y)^T\oplus(\lambda,\lambda)^T\in\cX_1$) 
if and only if $y\geq x+1$ (that is, $(x,y)^T\in \cX_1$), 
we get $\cX_1 \ominus \cB=\cX_1$. Therefore, 
\begin{eqnarray*}
A^{-1}(\cX_1 \ominus \cB) & = & A^{-1}(\cX_1) \\
& = & \set{(x,y)^T\in \zmax^2}{A(x,y)^T\in\cX_1} \\
& = & \set{(x,y)^T\in \zmax^2}{(y,x)^T\in\cX_1} \\
& = & \set{(x,y)^T\in \zmax^2}{x\geq y+1 }\;,
\end{eqnarray*}
and thus
\begin{eqnarray*}
{\cX}_2 & = & {\cX}_1 \cap A^{-1}({\cX}_1 \ominus {\cB})  \\
& = & \set{(x,y)^T\in \zmax^2}{y\geq x+1}  \cap  
\set{(x,y)^T\in \zmax^2}{x\geq y+1 }  \\
& = & \{ (-\infty ,-\infty )^T\} \; .
\end{eqnarray*}
Then, since by Lemma~\ref{lemaalgoAB} the sequence of semimodules 
$\{\cX_r\}_{r\in \N}$ is decreasing, it follows that 
$\cX_k=\{ (-\infty,-\infty)^T\}$ for all $k\geq 2$. Therefore, 
the maximal (geometrically) $(A,B)$-invariant semimodule contained 
in $\cK$ is trivial: $\cK^*=\cX_\omega=\{ (-\infty,-\infty)^T\} $. 
\end{example}

In the case of the theory of linear dynamical systems over a field, 
the sequence $\{\cX_r\}_{r\in \N}$ always converges in at most 
$n$ steps, since it is a decreasing sequence of subspaces 
of a vector space of dimension $n$. However, one of the problems in 
the max-plus case, which is reminiscent of difficulties of the theory 
of linear dynamical systems over rings 
(see~\cite{assan,AssLafPer,conte94,conte95,hautus82,hautus84}), 
is that the sequence $\{\cX_r\}_{r\in \N}$ may not stabilize 
(see Example~\ref{ejemplo2} below). This difficulty comes from 
the fact that the semimodule $\zmax^n$ is not Artinian, 
that is, there are infinite decreasing sequences of subsemimodules 
of $\zmax^n$. In the case of linear dynamical systems over rings, 
the convergence of the sequence $\{\cX_r\}_{r\in \N}$ in a 
finite number of steps is not guaranteed either, and 
although there exists a procedure for finding $\cK^*$ 
when $\cS$ is a Principal Ideal Domain (see~\cite{conte94}), 
in general the computation of $\cK^*$ remains a difficult problem.

\begin{example}\label{ejemplo2}
Let $\cS=\zmax$. Let us consider the matrices
\[
A=
\begin{pmatrix} 
-1 & -\infty \\ 
-\infty & 0 
\end{pmatrix}
\enspace \mbox{ and } \enspace
B=
\begin{pmatrix} 
0 \\ 
0
\end{pmatrix}\; ,
\]
and the semimodule $\cK=\set{(x,y)^T\in \zmax^2}{y\leq x-1}$. 
Note that $\cK=\im K$, where 
\[
K=
\begin{pmatrix} 
0 & 0 \\ 
-1 & -\infty 
\end{pmatrix}\;.
\] 
Next we show that in this case the sequence of semimodules 
$\{\cX_r\}_{r\in \N}$ defined by~\eqref{algoABinv} is given by: 
\begin{equation}\label{suceje2}
\cX_r=\set{(x,y)^T\in \zmax^2}{y\leq x-r}=
\im 
\begin{pmatrix} 
0 & 0 \\ 
-r & -\infty  
\end{pmatrix}
\;, 
\end{equation} 
for all $r\in \N$. We prove~\eqref{suceje2} 
by induction on $r$. Let us note, in the first place, 
that~\eqref{suceje2} is satisfied by definition when $r=1$. 
Assume now that~\eqref{suceje2} holds for $r=k$, that is: 
\[
\cX_k=\set{(x,y)^T\in \zmax^2}{y\leq x-k}=
\im 
\begin{pmatrix} 
0 & 0 \\ 
-k & -\infty  
\end{pmatrix}\;.
\] 
Let us note that $\cX_k \ominus \cB =\cX_k$, since 
there exists $\lambda \in \zmax$ such that 
$\max(y,\lambda)\leq\max(x,\lambda)-k$ (that is, 
there exists $(\lambda,\lambda)^T\in \cB$ such that 
$(x,y)^T\oplus(\lambda,\lambda)^T\in\cX_k$) 
if and only if $y\leq x-k$ (that is, $(x,y)^T\in \cX_k$). 
Therefore, 
\begin{eqnarray*}
A^{-1}(\cX_k \ominus \cB) & = & A^{-1}(\cX_k) \\
& = & \set{(x,y)^T\in \zmax^2}{A(x,y)^T\in\cX_k} \\
& = & \set{(x,y)^T\in \zmax^2}{(x-1,y)^T\in\cX_k} \\
& = & \set{(x,y)^T\in \zmax^2}{y\leq x-1-k }\;,
\end{eqnarray*}
and thus
\begin{eqnarray*}
\cX_{k+1} & = & \cX_k \cap A^{-1}(\cX_k \ominus \cB) \\
& = & \set{(x,y)^T\in \zmax^2}{y\leq x-k} \cap 
\set{(x,y)^T\in \zmax^2}{y\leq x-1-k } \\
& = & \set{(x,y)^T\in \zmax^2}{y\leq x-(1+k) }\;,
\end{eqnarray*}
which shows that~\eqref{suceje2} holds for all $r\in \N$. 

We see in this way that the sequence of semimodules 
$\{\cX_r\}_{r\in \N}$ is strictly decreasing and 
therefore does not stabilize. Let us finally note that the semimodule 
$\cX_{\omega}=\cap_{r\in \N} \cX_r=\set{(x,y)^T\in \zmax^2}{y=-\infty}$ 
is $A$-invariant, that is, $A(\cX_{\omega})\subset\cX_{\omega}$. 
Then, $\cX_{\omega}$ is in particular (geometrically) 
$(A,B)$-invariant and therefore   
$\cK^*= \cX_{\omega}=\set{(x,y)^T\in \zmax^2}{y=-\infty}$. 
\end{example}

An open problem is to determine whether it is always the case 
that $\cK^*= \cX_{\omega}$. 
It is worth mentioning that this equality does not necessarily hold in 
the case of linear dynamical systems over rings. 

\begin{remark}
Even when $\cS$ is a 
Principal Ideal Domain, it could be necessary to compute more than once 
(but a finite number of times) the limit $\cX_{\omega}$ of sequences 
defined as in~\eqref{algoABinv}. To be more precise, in such a case 
$\cX_1$ is defined as $\cK$ in the first step and, if it is necessary 
(that is, when $\cX_{\omega}$ is not a geometrically $(A,B)$-invariant module), 
in the next steps $\cX_1$ is defined as the smallest {\em closed} submodule 
containing the previous limit $\cX_{\omega}$ (see~\cite{conte94} for details). 
\end{remark}

Sufficient conditions for the stabilization of the 
sequence $\{\cX_r\}_{r\in \N}$ defined by~\eqref{algoABinv}, 
and therefore for the equality $\cK^*= \cX_{\omega}$ 
to hold true, will be given in Section~\ref{finitevolumeSec} 
in the case $\cS=\zmax$. Note that Example~\ref{ejemplo2} shows that 
even in the case of the tropical semiring $\nmin$ the sequence of 
semimodules $\{\cX_r\}_{r\in \N}$ may not stabilize (indeed all the 
computations in Example~\ref{ejemplo2} are valid when we restrict 
ourselves to the semiring $\N_{\max }^{-}=(\N^{-}\cup\{-\infty\},\max ,+)$, 
which is clearly isomorphic to $\nmin$). However, more general sufficient 
conditions for the equality $\cK^*= \cX_{\omega}$ to hold true can be 
given in the case of the tropical semiring using compactness arguments. 
With this aim, let us consider the topology of $\nmin $ defined by the 
metric:
\[
d(x,y)=|\exp(-x)-\exp(-y)| \; ,
\]
for all $x,y\in \nmin $. Note that $\nmin $ is compact equipped 
with this topology and therefore $\nmin^n$ is also compact 
equipped with the product topology. As a matter of fact, given 
a sequence $\{x_r\}_{r\in \N}$ of elements of $\nmin $, if the value $+\infty $ appears in 
$\{x_r\}_{r\in \N}$ an infinite number of times or if the set of finite values (that is, in $\N $) 
of $\{x_r\}_{r\in \N}$ is unbounded (in the usual sense), then $+\infty $ is an 
accumulation point of $\{x_r\}_{r\in \N}$. Otherwise, some finite element $x_k$ of 
$\{x_r\}_{r\in \N}$ must appear in this sequence an infinite number of times and then $x_k$ 
is an accumulation point of $\{x_r\}_{r\in \N}$. Now we have the following lemma.
\begin{lemma}\label{compact}
Finitely generated subsemimodules of $\nmin^n$ are compact.
\end{lemma}

\begin{proof}
Firstly, let us notice that $\nmin $ is a {\em topological semiring}, 
that is, for all sequences $\{x_r\}_{r\in \N}$ and $\{y_r\}_{r\in \N}$ 
of elements of $\nmin $ the following equalities are satisfied:
\[
\lim_{r\rightarrow \infty}\left(x_r\oplus y_r\right) = 
\left(\lim_{r\rightarrow \infty}x_r\right) \oplus 
\left(\lim_{r\rightarrow \infty}y_r\right) \; ,
\]
and
\[ 
\lim_{r\rightarrow \infty}\left(x_r\otimes y_r\right) = 
\left(\lim_{r\rightarrow \infty}x_r\right) \otimes 
\left(\lim_{r\rightarrow \infty}y_r\right) \; .
\] 
Let us now see that a finitely generated semimodule 
$\cX\subset \nmin^n$ is compact. Indeed, since $\cX$ is 
finitely generated there exists a matrix $Q\in \nmin^{n\times p}$, 
for some $p\in \N$, such that $\cX =\im Q$. Let $\{Qy_r\}_{r\in \N}$ 
be an arbitrary sequence of elements of $\cX$. To prove that $\cX$ is 
compact, we must show that $\{Qy_r\}_{r\in \N}$ has a subsequence 
which converges to an element of $\cX$. Since $\nmin^p$ is compact, we know that 
there exists a subsequence $\{y_{r_k}\}_{k\in \N}$ of $\{y_r\}_{r\in \N}$ 
and an element $y\in \nmin^p$ such that $\lim_{k\rightarrow \infty}y_{r_k}=y$. 
Then, using the fact that $\nmin$ is a topological semiring, it follows that 
\[
\lim_{k\rightarrow \infty}\left( Qy_{r_k}\right) = 
Q\left( \lim_{k\rightarrow \infty}y_{r_k}\right) = Qy\in \cX \; .
\]
Therefore, $\cX$ is compact.
\end{proof}

The following theorem shows that in the case of $\nmin$ the equality 
$\cK^*= \cX_{\omega}$ holds when $\cK$ is finitely generated.

\begin{theorem}
Let $\cK \subset \nmin^n$ be a finitely generated semimodule. Then, for all 
matrices $A\in \nmin^{n\times n}$ and $B\in \nmin^{n\times q}$, the maximal 
(geometrically) $(A,B)$-invariant semimodule $\cK^*$ contained in $\cK$ 
is given by $\cX_{\omega}=\cap_{r\in \N} \cX_r$, where the sequence 
of semimodules $\{\cX_r\}_{r\in \N}$ is defined by~\eqref{algoABinv}.
\end{theorem}
\begin{proof}
By Lemmas~\ref{obs1} and~\ref{lemaalgoAB}, to prove the theorem, it suffices to show 
that $\cX_{\omega}$ is a (geometrically) $(A,B)$-invariant semimodule, 
which is equivalent to showing that $\cX_{\omega}=\varphi(\cX_{\omega})$ 
by Lemma~\ref{obs2}. 

Since $\cX_{\omega}\subset \cX_r$ for all $r\in \N$, it follows that 
$\varphi(\cX_{\omega})\subset \varphi(\cX_r)=\cX_{r+1}$ for all 
$r\in \N$. Therefore, $\varphi(\cX_{\omega})\subset \cap_{r\in \N} \cX_r =
\cX_{\omega}$.

Let us now see that $\cX_{\omega}\subset \varphi(\cX_{\omega})$. 
Let $x$ be an arbitrary element of $\cX_{\omega}$. Then, since 
$x\in \varphi(\cX_r)=\cX_{r+1}$ for all $r\in \N$, we know that 
there exists a sequence $\{b_r\}_{r\in \N}\subset \cB$ such that 
$Ax\oplus b_r$ belongs to $\cX_r$ for all $r\in \N$. 
As $\cB$ is compact by Lemma~\ref{compact}, there exists $b\in \cB$ 
and a subsequence $\{b_{r_k}\}_{k\in \N}$ of $\{b_r\}_{r\in \N}$ such 
that $\lim_{k\rightarrow \infty}b_{r_k} =b$. Now, since by 
Lemma~\ref{lemaalgoAB} the sequence of semimodules 
$\{\cX_r\}_{r\in \N}$ is decreasing, 
it follows that $Ax\oplus b_{r_j}\in \cX_{r_k}$ for all $j\geq k$. 
Therefore, $Ax\oplus b\in \cX_{r_k}$ for all $k\in \N$ (recall that the 
semimodules $\cX_r$ are all finitely generated and then, 
by Lemma~\ref{compact}, in particular closed). Then, $Ax\oplus b$ 
belongs to $\cX_{\omega}$, from which we see that $x\in \varphi(\cX_{\omega})$. 
Therefore, $\cX_{\omega} \subset \varphi(\cX_{\omega})$.
\end{proof}

\section{Volume}\label{volumeSec}

In the next section we will give sufficient conditions on 
the semimodule $\cK$, when $\cS=\zmax$, to assure that the sequence 
of semimodules $\{\cX_r\}_{r\in \N}$ defined 
by~\eqref{algoABinv} stabilizes. For this purpose it is 
convenient to introduce first the notion of volume of a 
subsemimodule of $\zmax^n$ and study its properties.

\begin{definition}\label{defvolumen}
Let $\cK\subset \zmax^n$ be a semimodule. We call the {\rm volume} of 
$\cK$, represented by $\vol(\cK)$, the cardinality of the set 
$\set{x\in \cK}{x_1\oplus \cdots \oplus x_n=0}$, that is, 
$\vol(\cK)=\card\left(\set{x\in \cK}{x_1\oplus \cdots \oplus x_n=0}\right)$. 
Also, if $K\in \zmax^{n\times p}$, we represent by $\vol(K)$ 
the volume of the semimodule $\cK=\im K$, that is, 
$\vol(K)=\vol(\im K)$. 
\end{definition}

Before stating the following results, which provide some properties 
of the volume, it is convenient to introduce the following notation: 
if $\cX \subset \zmax^n$, then we define
$\tilde{\cX}=\set{x\in \cX}{x_1\oplus \cdots \oplus x_n=0}$.

\begin{remark}\label{obsvolproy}
Let us consider the max-plus parallelism relation 
$\sim$ on $\zmax^n$ defined by: $x\sim y$ if and only if 
$x=\lambda y$ for some $\lambda \in \R$ (that is, 
$x_i=\lambda+y_i$ for all $1\leq i\leq n$, in the usual algebra). 
We denote by $\cK/\sim$ the quotient of a semimodule $\cK\subset \zmax^n$ 
by this relation and by $[x]$ the equivalence class of $x\in \zmax^n$. 
Then, since the function $f:\tilde{\cK}\mapsto (\cK/\sim)-[\varepsilon]$ defined by 
$f(x)=[x]$ is a bijection, it follows that the volume of $\cK$ is equal to $\card(\cK/\sim)-1$, 
that is, the cardinality of the set of {\em nontrivial lines} (i.e. the equivalence classes of 
nonzero elements) contained in $\cK$. 
The {\em max-plus projective space} is the quotient of $\rmax^n$ by the parallelism relation.  
\end{remark}

\begin{lemma}\label{lemapropvol}
Let $A\in \zmax^{r\times n}$, $B\in \zmax^{n\times p}$ 
and $C\in \zmax^{p\times q}$ be matrices 
and $\cZ,\cY\subset  \zmax^n$ be semimodules. Then we have:
\begin{enumerate}
\item \label{p1} $\cY \subset \cZ \Rightarrow \vol(\cY) \leq \vol(\cZ) \;$,
\item \label{p2} if $\vol(\cY)<\infty$, 
then $\cY \varsubsetneq  \cZ \Rightarrow \vol(\cY) <\vol(\cZ) \;$,
\item \label{p3} $\vol(A\cY) \leq \vol(A)$ and then 
$\vol(AB) \leq \vol(A)\;$,
\item \label{p4} $\vol( A\cY) \leq \vol(\cY)$ and then 
$\vol(AB) \leq \vol(B)\;$,
\item \label{p5} $\vol(ABC) \leq \vol(B)\;$,
\item \label{p6} if $P\in \zmax^{n\times n}$ and $Q\in \zmax^{p\times p}$ 
are invertible\footnote{A matrix $P$ is invertible if there exists 
a matrix $P^{-1}$ such that $PP^{-1}=P^{-1}P=I$, where $I$ is the max-plus identity 
matrix. In the max-plus semiring, this means that the columns of $P$ are equal, 
up to a permutation, to the columns of $I$ multiplied by  non-zero scalars.}, 
then $\vol(PBQ) =\vol(B)\;$,
\item \label{p7} $\vol(A) =\vol( A^{T})\; $.
\end{enumerate}
\end{lemma}

\begin{proof}
\ref{p1}. This property is a consequence of the definition of volume: 
$\cY \subset \cZ \Rightarrow  \tilde{\cY} \subset \tilde{\cZ} 
\Rightarrow 
\card (\tilde{\cY})\leq \card (\tilde{\cZ}) \Rightarrow 
\vol(\cY)\leq  \vol(\cZ)$.

\ref{p2}. In the first place, we will show that the following simple 
property is satisfied: for all semimodules $\cY,\cZ\subset \zmax^n$, 
\begin{eqnarray}\label{tonta}
\cY \varsubsetneq \cZ \Rightarrow \tilde{\cY} \varsubsetneq \tilde{\cZ}\; .
\end{eqnarray}
As a matter of fact, assume that $\cY \varsubsetneq \cZ$. Then, there exists 
$x\in \cZ - \cY$. Therefore, we know that 
$x \neq(-\infty,\ldots ,-\infty)^T$ and we can define the 
vector $\tilde{x}=\left( x_1\oplus \cdots \oplus x_n\right)^{-1}x$ 
(that is, $\tilde{x}_i=x_i-\max\{x_1,\ldots ,x_n\}$ for all 
$1\leq i \leq n$, in the usual algebra). Now, it follows that 
$\tilde{x}\in \tilde{\cZ}- \tilde{\cY}$ and thus 
$\tilde{\cY} \varsubsetneq \tilde{\cZ}$. 
This proves property~\eqref{tonta}.

Now, using property~\eqref{tonta} and the fact that $\vol(\cY)<\infty$, 
we get: 
$\cY \varsubsetneq  \cZ \Rightarrow  
\tilde{\cY} \varsubsetneq \tilde{\cZ} \Rightarrow 
\card (\tilde{\cY}) < \card (\tilde{\cZ}) \Rightarrow  
\vol(\cY) < \vol(\cZ)$.

\ref{p3}. Since $A\cY\subset \im A$, applying Statement~\ref{p1}, 
we have: $\vol(A\cY) \leq \vol(\im A)=\vol(A)$.

\ref{p4}. From the definition of the set $\tilde{\cY}$ it follows 
that for each $y\in \cY- \{ (-\infty,\ldots ,-\infty )^T\}$ 
there exists $\tilde{y}\in \tilde{\cY}$ and $\lambda \in {\Z}$ such that 
$y=\lambda \tilde{y}$ (it suffices to take 
$\lambda = y_1\oplus \cdots \oplus y_n$ and $\tilde{y}=\lambda^{-1}y$). 
Therefore, 
\[
A\cY - \{ (-\infty,\ldots ,-\infty )^T\} \subset 
\set{\lambda A\tilde{y}}{\tilde{y}\in \tilde{\cY}, \lambda \in \Z}\;,
\]
and then we get:
\begin{eqnarray*}
&\vol(A\cY)=\card ( \set{x\in A\cY}{x_1\oplus \cdots \oplus x_r=0} ) \\
& \leq \card (\set{x=\lambda A\tilde{y}}{\tilde{y}\in \tilde{\cY}, \lambda \in \Z, x_1\oplus \cdots \oplus x_r=0} ) \\ 
& \leq \card (\set{A\tilde{y}}{\tilde{y}\in \tilde{\cY}}) \leq 
\card(\tilde{\cY})=\vol(\cY)\;.
\end{eqnarray*}

\ref{p5}. Applying Statements~\ref{p3} and~\ref{p4} we get: 
$\vol(ABC)\leq \vol(AB)\leq \vol(B)$.

\ref{p6}. From Statement~\ref{p5} we obtain: 
$\vol(B) =\vol(P^{-1}PBQQ^{-1})\leq \vol(PBQ)\leq \vol(B)$. 
Therefore, $\vol(B) =\vol(PBQ)$.

\ref{p7}. Let us note, in the first place, 
that we can define in a completely analogous way 
the volume of a subsemimodule of $\zmin^n$. Then, 
since the function $x\rightarrow -x$ is an 
isomorphism from $\zmax$ to $\zmin$, it is clear that 
$\vol(\cZ)=\vol(-\cZ)$ for every subsemimodule 
$\cZ\subset \zmax^n$. Let us now consider the matrix 
$A^{\sharp}=-A^T$ and the semimodule 
$\cY=\im(A^{\sharp})\subset \zmin^n$. Since 
$\cY=-\im(A^T)$, we know that $\vol(A^T)=\vol(\cY)$. 
Now, using elements of residuation theory (we refer the reader to~\cite{BlythJan72} 
for an extensive presentation of this theory), it can be shown 
(see for example~\cite{bcoq} or~\cite{gaubert01a}) 
that the following two properties hold:
\begin{eqnarray*}
A(A^{\sharp}(Ax)) & = & Ax\;,\enspace \forall x\in \zmax^n\;,  \mbox{ and }\\
A^{\sharp}(A(A^{\sharp}y)) & = & A^{\sharp}y\;,\enspace \forall y\in \zmin^r \;,
\end{eqnarray*}
where the products by $A$ are performed in $\zmaxb$ and the products 
by $A^{\sharp}$ are performed in $\zminb$. Therefore, 
the function $f:\im(A)\mapsto \im(A^{\sharp})$ defined by 
$f(y)=A^{\sharp}y$ is a bijection with inverse $g(x)=Ax$. 
Then, the function $F$ from $\im(A)/\sim$ to $\im(A^{\sharp})/\sim$ 
defined by $F([y])=[A^{\sharp}y]$, where $[x]$ denotes 
the equivalence class of $x$ by the parallelism relation $\sim$, 
is also a bijection. Now, using Remark~\ref{obsvolproy}, 
we obtain: 
$\vol(A)=\card(\im(A)/\sim)-1=\card(\im(A^{\sharp})/\sim)-1=\vol(A^{\sharp})=\vol(\cY)$, 
and then $\vol(A)=\vol(\cY)=\vol(A^T)$.
\end{proof}

\section{Specifications with finite volume}\label{finitevolumeSec}

In the next theorem we give a condition on the specification $\cK$, 
when $\cS=\zmax$, ensuring that the sequence of semimodules defined 
by~\eqref{algoABinv} stabilizes. 

\begin{theorem}\label{th-inv}
Let $\cK\subset\zmax^n$ be a semimodule with finite volume. 
Then, for all $A\in\zmax^{n\times n}$ and $B\in \zmax^{n\times p}$, 
the maximal (geometrically) $(A,B)$-invariant semimodule $\cK^*$ 
contained in $\cK$ is finitely generated. Moreover, 
if we define the sequence of semimodules $\{\cX_r\}_{r\in \N}$ 
by~\eqref{algoABinv}, then $\cK^*=\cX_k$ for some $k\leq \vol(\cK)+1$.
\end{theorem}

\begin{proof}
First of all, let us note that every semimodule $\cY\subset\zmax^n$ 
with finite volume is necessarily finitely generated. Indeed, 
this property is a consequence of the fact that 
$\cY=\smallvect(\tilde{\cY})$. Now, as $\cK^*\subset\cK$, applying 
Statement~\ref{p1} of Lemma~\ref{lemapropvol} it follows that 
$\vol(\cK^*)\leq \vol(\cK)<\infty$, and then $\cK^*$ is finitely generated. 

Let us now see that the sequence of semimodules 
$\{\cX_r\}_{r\in \N}$ defined by~\eqref{algoABinv} 
must stabilize in at most $\vol(\cK)+1$ steps. 
Indeed, by Lemma~\ref{lemaalgoAB} 
we know that the sequence of semimodules 
$\{\cX_r\}_{r\in \N}$ is decreasing. Then, 
using Statement~\ref{p1} of Lemma~\ref{lemapropvol}, 
we see that $\{\vol(\cX_r)\}_{r\in \N}$ is a decreasing 
sequence of nonnegative integers. Therefore, 
there exists $k\leq \vol(\cX_1)+1=\vol(\cK)+1$ such that 
$\vol(\cX_{k+1})=\vol(\cX_k)$. Then, 
as $\cX_{k+1}\subset \cX_k\subset \cK$ by Lemma~\ref{lemaalgoAB}, 
we know that $\vol(\cX_{k+1})=\vol(\cX_k)\leq \vol(\cK)<\infty$ 
(once again, by Statement~\ref{p1} of Lemma~\ref{lemapropvol}). 
Finally, applying Statement~\ref{p2} of Lemma~\ref{lemapropvol} 
to the semimodules $\cX_{k+1}$ and $\cX_k$, it follows that  
$\cX_{k+1}=\cX_k$, from which we conclude that $\cK^*=\cX_k$.
\end{proof}
   
An important particular case of Theorem~\ref{th-inv} 
is the one in which the semimodule $\cK$ is 
generated by a finite number of vectors whose 
entries are all finite. In this case it is possible 
to bound the volume of $\cK$ by means of the additive 
version of Hilbert's projective metric: 
for all $x\in\Z^n$, define 
\[
\|x\|_H=\max\set{x_i}{1\leq i \leq n}-\min\set{x_i}{1\leq i \leq n} \enspace,
\]
and for all $K\in \Z^{n\times s}$, define 
\[
\Delta_H(K)= \max\set{\|K_{\cdot i}\|_H}{1\leq i\leq s} \enspace,
\]
where $K_{\cdot i}$ denotes the $i$-th column of the matrix $K$. 
Then we have the following corollary.

\begin{corollary}\label{corvolfin} 
Let $\cK=\im K$, where $K\in \zmax^{n\times s}$ 
is a matrix whose entries are all finite. 
Then, for all $A\in\zmax^{n\times n}$ and $B\in \zmax^{n\times p}$, 
the maximal (geometrically) $(A,B)$-invariant semimodule $\cK^*$ 
contained in $\cK$ is finitely generated and, if we define the 
sequence of semimodules $\{\cX_r\}_{r\in \N}$ by~\eqref{algoABinv}, 
there exists some $k\leq (\Delta_H(K)+1)^n-\Delta_H(K)^n+1$ 
such that $\cK^*=\cX_k$.
\end{corollary}

\begin{proof}
By Theorem~\ref{th-inv}, to prove the corollary, 
it suffices to show that 
\begin{equation}\label{projective}
\vol(\cK)\leq (\Delta_H(K)+1)^n-\Delta_H(K)^n  \enspace ,
\end{equation}
where the power $n$ is in the usual algebra.

Since the additive version of Hilbert's projective metric 
$\|\cdot\|_H$ satisfies the following properties:
\begin{eqnarray*}
\|\lambda x\|_H & = & \|x\|_H \;, \\
\|x\oplus y\|_H & \leq &\|x\|_H\oplus \|y\|_H \;,
\end{eqnarray*}
for all $x,y\in \Z^n$ and $\lambda \in \Z$, 
it follows that $\|x\|_H\leq\Delta_H(K)$ for all 
$x\in \cK- \{(-\infty,\ldots ,-\infty)^T\}$ 
and therefore $\cK$ is contained in the semimodule 
\[
\cY=\set{x\in \Z^n}{\|x\|_H\leq\Delta_H(K)}\cup \{(-\infty,\ldots ,-\infty)^T\} 
\] 
(note that the only vector in $\cK$ with at least one entry 
equal to $-\infty$ is $(-\infty,\ldots ,-\infty)^T$). 
Then, by Statement~\ref{p1} of Lemma~\ref{lemapropvol}, 
to prove~\eqref{projective} it suffices 
to show that $\vol(\cY)=(\Delta_H(K)+1)^n-\Delta_H(K)^n$. 
With this aim, we must compute the number of elements of the set:
\begin{eqnarray*}
\tilde{\cY} & = & \set{x\in \cY}{x_1\oplus \cdots \oplus x_n=0}  \\
& = & \set{x\in \Z^n}{\|x\|_H\leq\Delta_H(K),x_1\oplus \cdots \oplus x_n=0}\;,
\end{eqnarray*}
that is, the number of vectors $x$ in $\Z^n$ with entries 
between $-\Delta_H(K)$ and zero (since $\max_ix_i=x_1\oplus \cdots \oplus x_n=0$ and 
$\Delta_H(K)\geq \|x\|_H =\max_ix_i-\min_ix_i=-\min_ix_i$) and with at least one entry 
equal to zero (since $\max_ix_i=0$). 
We know that there are ${n \choose r} \Delta_H(K)^{n-r}$ 
elements in the set $\tilde{\cY}$ with exactly $r$ entries equal to zero. 
To be more precise, there exist ${n \choose r}$ 
different ways of choosing the $r$ entries which will 
have the value zero, and there exist $\Delta_H(K)^{n-r}$ 
different ways of assigning values to the $n-r$ 
remaining entries among the $\Delta_H(K)$ possible values. 
Therefore, the number of elements of the set $\tilde{\cY}$ is:
\[
\sum_{r=1}^{r=n} {n \choose r}\Delta_H(K)^{n-r} = 
(\Delta_H(K)+1)^n-\Delta_H(K)^n \;, 
\] 
and then $\vol(\cY)=(\Delta_H(K)+1)^n-\Delta_H(K)^n$.
\end{proof}

Note that in the proof of Corollary~\ref{corvolfin} we showed, 
in particular, that for each matrix $K\in \zmax^{n\times s}$ 
whose entries are all finite, the volume $\vol(K)$ is bounded 
by $(\Delta_H(K)+1)^n-\Delta_H(K)^n$ (this is inequality~\eqref{projective}). 
We next show that this bound is tight. Indeed, 
let us consider the semimodule 
\[
\cY=\set{x\in \Z^n}{\|x\|_H\leq M}\cup \{(-\infty,\ldots ,-\infty)^T\}\;,
\]
where $M\in \N$. Note that in the proof of Corollary~\ref{corvolfin} 
we proved that $\cY$ has volume $(M+1)^n-M^n$. Now, 
if we define the matrix $K\in \zmax^{n\times n}$ by 
$K_{ij}=M$ if $i=j$ and $K_{ij}=0$ otherwise, 
it follows that $\cY=\im(K)$ and $\Delta_H(K)=M$. 
Therefore, there exist matrices  $K\in \zmax^{n\times s}$ 
(whose entries are all finite) which have volume equal to 
$(\Delta_H(K)+1)^n-\Delta_H(K)^n$.

Theorem~\ref{th-inv} is useful in many practical problems because 
in such problems the specification $\cK$ frequently has finite volume. 
This is often the case when $\cK$ models certain stability conditions, 
as for example, ``bounded delay'' requirements. To be more precise, let us 
assume that system~\eqref{dynamicsystem} is the dater representation 
of a timed event graph (we refer the reader to~\cite{bcoq} for more 
details on the modeling of timed event graphs). Then, a typical 
case of semimodule $\cK$ which arises in applications is: 
\begin{equation}\label{semiacot}
\cK =\set{x\in \zmax^n}{x_i-x_j \leq d_{ij}, \forall 1\leq i,j\leq n}\; ,
\end{equation}
where $D=(d_{ij})$ is a matrix with entries in $\Z\cup \{+\infty\}$. 
Note that the state vector $x(k)$, representing the dates of the 
firings numbered $k$, belongs to $\cK$ if and only if $x(k)_i-x(k)_j \leq d_{ij}$, 
for all $1\leq i,j\leq n$, which means that the delay between the 
$k$-th firing of the transition labeled $j$ and the $k$-th 
firing of the transition labeled $i$ should not exceed $d_{ij}$. 
Note also that in practice we usually can assume that $D$ only has finite 
entries, since we can replace $+\infty$ by a sufficiently large constant. 
We next show that in such a case, the semimodule $\cK$ defined by~\eqref{semiacot} 
has finite volume. Let us first recall that a 
directed graph $\cG (A)$, called the {\em precedence graph} of $A$, 
is associated with a matrix $A=(a_{ij})\in \rmax^{n\times n}$. This graph 
is defined as follows: there exists a directed arc of 
{\em weight} $a_{ji}$ from node $i$ to node $j$ 
if and only if $a_{ji}\not = -\infty$. A matrix whose 
precedence graph is strongly connected is called 
{\em irreducible}. The spectral radius $\rho_{\max }(A)$ 
of $A$ is defined by:
\[
\rho_{\max }(A)=\bigoplus_{k=1}^{n}\mbox{tr}(A^k)^\frac{1}{k}= 
\max_{1\leq k\leq n} \max_{i_1,\ldots ,i_k} \frac{a_{i_1i_2}+\cdots 
+a_{i_ki_1}}{k}  \; ,
\]
that is, the maximal circuit mean of $\cG (A)$.

Before stating the following lemma, which shows in particular that the 
semimodule~\eqref{semiacot} has finite volume
when $D$ only has finite entries, let us note that    
\begin{equation}\label{semiacot2}
\cK =\set{x\in \zmax^n}{Ex \leq x}\; ,
\end{equation}
where $E=(-D)^T$. Then we have:

\begin{lemma}\label{lemaHstar}
If the matrix $E$ is irreducible, then the semimodule 
$\cK$ defined by~\eqref{semiacot2} has finite volume. 
Moreover, if $E$ has spectral radius strictly greater 
than the unit (that is, 0), then $\cK$ reduces to the 
null vector.
\end{lemma}

\begin{proof}
In the first place, let us see that 
$\cK=\im(E^*)\cap \zmax^n$, where 
\[
E^*=\bigoplus_{r=0}^{\infty }E^r=I\oplus E\oplus E^2\oplus \cdots 
\] 
(note that the matrix $E^*$ can have entries equal to $+\infty$, 
so that $E^*$ should be thought of as a map from $\zmaxb^n$ to $\zmaxb^n$).
Indeed, we have:
\begin{eqnarray*}
& x\in \cK  \Rightarrow  Ex\leq x,x\in \zmax^n \Rightarrow \\
& E^{r}x\leq x , \forall r\in \N , x\in \zmax^n 
\Rightarrow E^*x\leq x , x\in \zmax^n \Rightarrow  \\
&  E^*x = x , x\in \zmax^n \Rightarrow x\in \im(E^*)\cap \zmax^n \; ,
\end{eqnarray*}
and 
\begin{eqnarray*}
& x\in \im(E^*)\cap \zmax^n  \Rightarrow \\
& x=E^*y, \mbox { for some } y\in \zmaxb^n , x\in \zmax^n \Rightarrow \\
& Ex\leq E^*x=E^*E^*y=E^*y=x , x\in \zmax^n \Rightarrow x\in \cK \; .
\end{eqnarray*}

When $E$ has spectral radius less than or equal to the unit, we know that: 
\[
E^*=I\oplus E\oplus \cdots \oplus E^{n-1}\; ,
\]
since $E^r\leq I\oplus E\oplus \cdots \oplus E^{n-1}$ 
for all $r\geq n$ (see for example Theorem~3.20 of~\cite{bcoq}). 
Moreover, since $E$ is irreducible, we know that all the entries 
of $E^*$ are finite. Indeed, this follows from the fact that $E^k_{ij}$, 
for $i\not =j$, is the maximal weight of all paths of length $k$ running from 
$j$ to $i$ in the precedence graph of $E$. Then, 
the proof of Corollary~\ref{corvolfin} shows that $\cK$ has finite volume.

When $E$ has spectral radius strictly greater than the unit, since 
$E$ is irreducible, all the entries of $E^*$ are equal to 
$+\infty$ (once again by the interpretation of the entries of 
the matrix $E^k$ in terms of the weight of paths in the precedence 
graph of $E$). Therefore, the only vector in $\cK=\im (E^*)\cap \zmax^n$ 
is the null vector.
\end{proof}

We end this section with an example showing that in Theorem~\ref{th-inv}, 
the bound $\vol(\cK)+1$ on the number of steps needed to stabilize 
the sequence of semimodules $\{\cX_r\}_{r\in \N}$ defined by~\eqref{algoABinv}, 
cannot be improved. 
 
\begin{example}
Let us consider the matrices
\[
A=
\begin{pmatrix} 1 & -\infty \\ -\infty & 0 \end{pmatrix}
\enspace \mbox{ and } \enspace
B=\begin{pmatrix} 0 \\ 0\end{pmatrix}\; ,
\]
and the semimodule $\cK=\set{(x,y)^T\in \zmax^2}{x+1\leq y\leq x+l}$, 
where $l\in \N$. Then, in this case we have:
\[
\tilde{\cK}=\set{(x,y)^T\in \cK}{x\oplus y=0}=\{(-1,0)^T,\ldots ,(-l,0)^T\}\; , 
\]
from which we get $\vol(\cK)=l$. 
Therefore, we are able to apply Theorem~\ref{th-inv}. 
In fact, $\cK=\im K$ where 
\[
K=\begin{pmatrix} 0 & 0 \\ 1 & l \end{pmatrix}\;,
\]
so we are also in a position to apply Corollary~\ref{corvolfin}. 

By Theorem~\ref{th-inv} we know that the sequence of semimodules 
$\{\cX_r\}_{r\in \N}$ defined by~\eqref{algoABinv} must stabilize 
in at most $\vol(\cK)+1=l+1$ steps. 
Let us check this fact in this particular case. 
In the first place, note that  
$\cK\subset \set{(x,y)^T\in \zmax^2}{x+1\leq y}$, so that 
$\cX_r\subset\cK\subset \set{(x,y)^T\in \zmax^2}{x+1\leq y}$ 
for all $r\in \N$. 
Then, it is easy to show (applying a straightforward variant 
of the computation of $\cX_r\ominus \cB$ done in 
Example~\ref{ejemplo2}) that $\cX_r\ominus \cB=\cX_r$ 
for all $r\in \N$. In this way we get:
\begin{eqnarray*}
\cX_1 & = & \left\{(x,y)^T\in \zmax^2\mid x+1\leq y\leq x+l\right\} 
\;,\\
\cX_2 & = &\cX_1 \cap A^{-1}(\cX_1 \ominus \cB)=\cX_1 \cap A^{-1}(\cX_1) \\
& = &\left\{(x,y)^T\in \zmax^2\mid x+1\leq y\leq x+l\right\}\cap 
\\ & & \enspace \;
\left\{(x,y)^T\in \zmax^2\mid x+2\leq y\leq x+l+1\right\} \\
& = &\left\{(x,y)^T\in \zmax^2\mid x+2\leq y\leq x+l\right\} 
\varsubsetneq \cX_1\;,\\
& \vdots & \\
\cX_l & = & \cX_{l-1} \cap A^{-1}(\cX_{l-1} \ominus \cB)=\cX_{l-1} \cap A^{-1}(\cX_{l-1}) \\
& = & \left\{(x,y)^T\in \zmax^2\mid x+l-1\leq y\leq x+l\right\}\cap 
 \\ & & \enspace \;
\left\{(x,y)^T\in \zmax^2\mid x+l\leq y\leq x+l+1\right\} \\
& = & \left\{(x,y)^T\in \zmax^2\mid x+l\leq y\leq x+l\right\}  \\
& = & \left\{(x,y)^T\in \zmax^2\mid y= x+l\right\} 
\varsubsetneq \cX_{l-1} \;,\\
\cX_{l+1} & = &\cX_l \cap A^{-1}(\cX_l \ominus \cB)=\cX_l \cap A^{-1}(\cX_l)  \\
& = & \left\{(x,y)^T\in \zmax^2\mid y = x+l\right\}\cap 
\left\{(x,y)^T\in \zmax^2\mid y= x+l+1\right\}  \\
& = & \left\{(-\infty ,-\infty)^T\right\} \varsubsetneq \cX_l\;. 
\end{eqnarray*}
Then, since by Lemma~\ref{lemaalgoAB} we know that   
\[
\left\{(-\infty ,-\infty)^T\right\} \subset \cX_{l+2} \subset
\cX_{l+1}=\left\{(-\infty ,-\infty)^T\right\}\; , 
\]
it is clear that $\cX_{l+2}=\cX_{l+1}$, and therefore 
\[
\cK^*=\cX_{l+1}=\left\{(-\infty ,-\infty)^T\right\} \; .
\] 
In this way we see that in this particular case the sequence 
of semimodules $\{\cX_r\}_{r\in \N}$ stabilizes in exactly 
$\vol(\cK)+1=l+1$ steps.
\end{example}

\section{Algebraically $(A,B)$-invariant semimodules}\label{algABinvSec}

This section deals with another fundamental problem in the geometric 
approach to the theory of linear dynamical systems: 
the computation of a linear feedback. Let us once again consider the dynamical 
system~\eqref{dynamicsystem}. Let us assume that we already know the maximal 
(geometrically) $(A,B)$-invariant semimodule $\cK^*$ contained in a given 
semimodule $\cK\subset \cS^n$. From a dynamical point of view, this means 
that the trajectories of system~\eqref{dynamicsystem} starting in $\cK^*$ 
can be kept inside $\cK^*$ by a suitable choice of the control. Our new problem 
is to determine whether this control can be generated by using a state feedback. 
In other words, we want to determine whether there exists a linear feedback 
$u(k)=Fx(k-1)$, where $F\in \cS^{q\times n}$, which makes $\cK^*$ invariant 
with respect to the resulting closed loop system:
\begin{equation}\label{sisaut}
x(k)=(A\oplus BF)x(k-1)\; ,
\end{equation}
that is, such that every trajectory of the closed loop 
system~\eqref{sisaut} is completely contained in $\cK^*$ when its initial 
state is in $\cK^*$. If a linear feedback with this property exists, we will 
say that $\cK^*$ is an algebraically $(A,B)$-invariant semimodule. Some 
authors call this notion $(A+BF)$-invariance (see~\cite{assan}) or the 
feedback property (see~\cite{hautus82,conte95,conte94}).
\begin{definition}\label{defABFinv}
Given the matrices $A\in\cS^{n\times n}$ and $B\in \cS^{n\times q}$, we say 
that a semimodule $\cX \subset \cS^n$ is {\rm algebraically 
$(A,B)$-invariant} if there exists $F\in \cS^{q\times n}$ such that 
\[
(A\oplus BF) \cX \subset \cX \enspace .
\]
\end{definition}

Obviously, every algebraically $(A,B)$-invariant semimodule is 
also geometrically $(A,B)$-invariant. Nevertheless, when $\cS=\zmax$ 
it is not clear whether a geometrically $(A,B)$-invariant semimodule 
is algebraically $(A,B)$-invariant. 
Once again, this problem is reminiscent of difficulties of 
the theory of linear dynamical systems over rings 
(see~\cite{hautus82,hautus84,conte94,conte95,assan,AssLafPer}). 
Indeed, in the case of linear dynamical systems with 
coefficients in a field, the class of geometrically 
$(A,B)$-invariant spaces coincides with the class of 
algebraically $(A,B)$-invariant spaces (see~\cite{wonham}). 
This property makes the (geometrically) $(A,B)$-invariant 
spaces very useful in the classical theory. 
However, this crucial feature is no longer true for 
linear dynamical systems with coefficients in a ring, 
that is, there exist geometrically $(A,B)$-invariant 
modules which are not algebraically $(A,B)$-invariant 
(see~\cite{hautus82}, in particular Example~2.3). 
The following example shows that this is also the case for 
linear dynamical systems over the tropical 
semiring $\nmin=(\Np,\min,+)$.  

\begin{remark}
In the case of rings, a necessary and sufficient condition for $\cK^*$ 
to be algebraically $(A,B)$-invariant can be given in the form of a factorization 
condition on the transfer function, assuming that the system is 
reachable and injective (see~\cite{hautus82}). When $\cS$ is a 
Principal Ideal Domain, it can be shown that $\cK^*$ is 
algebraically $(A,B)$-invariant if and only if it is a direct 
summand (see~\cite{hautus82,conte95,conte94}).
\end{remark}

\begin{example}
Let $\cS=\nmin$. Let us consider the matrices
\[
A=
\begin{pmatrix} 1 & +\infty  \\ 1 & 0 \end{pmatrix}
\enspace \mbox{ and } \enspace
B=\begin{pmatrix} 1 \\ 1\end{pmatrix}\; ,
\]
and the semimodule $\cK=\left\{(x,y)^T\in \nmin^2\mid x\leq y \right\}$.

In the first place, let us compute the maximal geometrically 
$(A,B)$-invariant semimodule $\cK^*$ contained in $\cK$. With this aim, 
we will compute the sequence of semimodules $\{\cX_r\}_{r\in \N}$ 
defined by~\eqref{algoABinv}. We have:
\begin{eqnarray*}
\cX_1 & = & \cK=\left\{(x,y)^T\in \nmin^2\mid x\leq y \right\}
\;,\\
\cX_2 & = & \cX_1 \cap A^{-1}(\cX_1 \ominus \cB) \\ 
& = & \left\{(x,y)^T\in \nmin^2\mid x\leq y \right\}\cap 
\left\{(x,y)^T\in \nmin^2\mid 1\leq y\right\} \\
& = & \left\{(x,y)^T\in \nmin^2\mid x\leq y, 1\leq y\right\} 
\;, \\
\cX_3 & = & \cX_2 \cap A^{-1}(\cX_2 \ominus \cB)= \\
& = & \left\{(x,y)^T\in \nmin^2\mid x\leq y,1\leq y\right\}\cap 
\left\{(x,y)^T\in \nmin^2\mid 1\leq y\right\} \\
& = & \cX_2 \;. 
\end{eqnarray*}
Then, we get $\cK^*=\cX_2=\left\{(x,y)^T\in \nmin^2\mid x\leq y,1\leq y\right\}$. 
Indeed, it is easy to check that a trajectory which starts at a point of 
$\cK^*=\cK- \left\{(0,0)^T\right\}$ can be kept inside $\cK$ with 
the sequence of controls identically equal to $(1,1)^T$, and that a 
trajectory which starts at the point $(0,0)^T$ cannot be kept inside $\cK$ 
(since for all controls in $\cB$ the next state of the system is always 
$(1,0)^T$, which does not belong to $\cK$).

Let us now see that $\cK^*$ is not an algebraically $(A,B)$-invariant 
semimodule. With this aim, we will show that a trajectory which starts at the 
point $(1,1)^T\in\cK^*$ cannot be kept inside $\cK^*$ when a linear state 
feedback is applied. Let $F\in \nmin^{1\times2}$ be an arbitrary feedback. 
Then, since $F (1,1)^T\geq 1$, we know that $BF(1,1)^T=(\alpha,\alpha)^T$, 
where $\alpha \geq 2$. Therefore, 
$$(A\oplus BF)\begin{pmatrix} 1 \\ 1\end{pmatrix} = 
\begin{pmatrix} 2 \\ 1\end{pmatrix}\oplus \begin{pmatrix} \alpha \\ \alpha \end{pmatrix}=
\begin{pmatrix} 2 \\ 1\end{pmatrix}\not \in \cK^*\; ,$$
which shows that $\cK^*$ is not an algebraically $(A,B)$-invariant semimodule.
\end{example}

We next show how we can decide, using the existing results 
on max-plus linear system of equations, whether a finitely generated 
subsemimodule of $\zmax^n$ is algebraically $(A,B)$-invariant. 
This method also computes a linear feedback with the required 
property when it exists. Let $A\in\zmax^{n\times n}$, $B\in \zmax^{n\times q}$, 
and let $\cX$ be a finitely generated subsemimodule of 
$\zmax^n$, so that there exists $Q\in \zmax^{n\times r}$, 
for some $r\in \N$, such that $\cX=\im Q$. Then, 
from Definition~\ref{defABFinv} it readily follows 
that $\cX$ is an algebraically $(A,B)$-invariant semimodule 
if and only if there exist matrices 
$F\in \zmax^{q\times n}$ and $G\in \zmax^{r\times r}$ 
such that:
\begin{equation}\label{sistnohom}
(A\oplus BF)Q=QG\; .
\end{equation}
As~\eqref{sistnohom} is a two sided max-plus linear system of equations, 
we know that its set of solutions $(F,G)$ is a finitely 
generated max-plus convex set, 
which can be explicitly computed by the general elimination methods 
(see~\cite{butkovicH,gaubert92a,gaubert98n,maxplus97}). 
In this way we see that we can effectively decide whether a finitely 
generated subsemimodule of $\zmax^n$ is algebraically $(A,B)$-invariant. 

\begin{remark}\label{obssisequa}
The elimination algorithm shows that the set of solutions of a homogeneous 
max-plus linear system of the form $Dx=Cx$, where $D, C$ are matrices of 
suitable dimensions, is a finitely generated semimodule. This algorithm relies 
on the fact that hyperplanes of $\rmax^n$ (that is, the set of solutions of
an equation of the form $dx=cx$, where $d,c\in \rmax^n$ are row vectors) 
are finitely generated. It is worth mentioning that the resulting 
naive algorithm has an a priori doubly exponential complexity. However, the 
doubly exponential bound is pessimistic. It is possible to incorporate 
in this algorithm the elimination of redundant generators which reduces its 
execution time. In fact, we are currently working on this subject and we 
believe that improvements are possible, since we have shown by direct 
arguments that the number of generators of the set of solutions is at 
most simply exponential. This will be the subject of a further work.
\end{remark}

Let us note that to decide whether $\cX=\im Q$ is an algebraically 
$(A,B)$-invariant semimodule it suffices to know whether the system of 
equations~\eqref{sistnohom} has at least one solution. Taking this into account, 
it is worth mentioning that there are algorithms to compute a single solution 
(with finite entries) of homogeneous max-plus linear systems which seem to be 
more efficient in practice than the elimination methods (see~\cite{bcg99,walkup}).  
Indeed, it is known that the problem of the existence of a solution (with finite entries) 
of a homogeneous max-plus linear system can be reduced to the problem of the existence 
of a sub-fixed point of a min-max function (for more background on min-max functions 
we refer the reader to~\cite{cras,gg} and the references therein). To be more precise, 
observe that $Dx=Cx$ is equivalent to 
$x\leq \min \left\{ D\backslash Cx, C\backslash Dx\right\}$, where 
$D\backslash Cx=\sup \set{y\in \rmaxb^n}{Dy\leq Cx}$ ($C\backslash Dx$ 
is defined analogously). Since $D\backslash Cx$ can be computed as 
$(-D^T)(Cx)$, where the product by $-D^T$ is performed in $\rminb$ (see~\cite{bcoq}), 
it follows that $f(x)=\min \left\{ D\backslash Cx, C\backslash Dx\right\}$ is a min-max 
function. Then, there is $x\in \R^n$ such that $x\leq f(x)$ 
(that is, a sub-fixed point of $f$) if and only if all the entries of the 
{\em cycle time vector} of $f$, which is defined as 
$\chi(f)=\lim_{k\rightarrow \infty} f^k(x)/k$, are nonnegative 
(see~\cite{cras,gg}). The cycle time vector $\chi(f)$, and, if it exists, a solution of 
$x\leq f(x)$ can be efficiently computed via the min-max Howard algorithm (we refer the 
reader to~\cite{cras,gg} for a detailed presentation of this algorithm). 
Although the min-max Howard algorithm behaves remarkably well in practice, 
its complexity is not yet well understood (\cite{cras,gg}). 

To be able to apply this algorithm to solve our problem, firstly we need to 
add one unknown $t$ to system~\eqref{sistnohom} in order to obtain a homogeneous max-plus 
linear system of equations:
\begin{equation}\label{sisthom}
(At\oplus BF)Q=QG\; .
\end{equation}
Then, as system~\eqref{sistnohom} has at least one solution if and only if 
system~\eqref{sisthom} has at least one solution with $t\not = -\infty$, 
the semimodule $\cX=\im Q$ is algebraically $(A,B)$-invariant if and only if 
system~\eqref{sisthom} has at least one solution with $t\not = -\infty$ 
(note that if $(t,F,G)$ is a solution of~\eqref{sisthom} with $t\not = -\infty$, 
then $t^{-1}F=(-t)F$ is the feedback we are looking for). 
Therefore, as $(t,F,G)$ is a solution of~\eqref{sisthom} if and only if
\begin{eqnarray}\label{subpunfijo} 
t & \leq & (AQ)\backslash (QG)\; , \nonumber \\ 
F & \leq & B\backslash (QG)/Q \; , \label{sistdes}\\
G & \leq & Q\backslash ((At\oplus BF)Q) \; , \nonumber
\end{eqnarray}
where $D\backslash C$ is defined as $\sup \set{E\in \zmaxb^{p\times r}}{DE\leq C}$ 
for all $D\in \zmax^{n\times p}$ and $C\in \zmax^{n\times r}$ (the function $/$ is 
defined in an analogous way), if we can find a sub-fixed point of the min-max function 
defined by the right hand side of~\eqref{sistdes}, then the semimodule $\cX=\im Q$ is 
algebraically $(A,B)$-invariant. 
 
\section{Application to transportation networks with a timetable}\label{aplicacionSec}

Let us consider the railway network given in Figure~\ref{figure1}. 
Firstly, we will recall how the evolution of this kind of 
transportation network can be described by max-plus linear 
dynamical systems of the form of~\eqref{dynamicsystem} (we refer the reader 
to~\cite{bcoq,OlsSubGett98,braker91,deVDeSdeM98} for details on max-plus models for 
transportation networks). We are interested in the departure times of the trains from the stations. 
Let us assume that in the initial state there is a train running along 
each of the following tracks: the one connecting $P$ with $Q$, the one connecting $Q$ with $P$, 
the one connecting $Q$ with $Q$ via $R$, and finally the one connecting $Q$ with $Q$ via $S$. 
We call these tracks directions $d_1$, $d_2$, $d_3$ and $d_4$ respectively, as it is shown in 
Figure~\ref{figure1}. In general, we can have $n$ different directions. The traveling time in 
direction $d_i$ (to which the time needed for passengers to leave and board the train is added) 
will be denoted by $t_i$. For our example these times are given in Figure~\ref{figure1}. 
Let $x_i(k)$ denote the $k$-th departure time of the train which leaves in direction $d_i$. 
As we explained in the introduction, a train cannot leave before a number of 
conditions have been satisfied. 
A first condition is that the train must have arrived at the station. 
For instance, let us assume that the train which leaves in 
direction $d_i$ is the one which comes from direction $d_{r(i)}$ 
(in Figure~\ref{figure1} we have: $r(1)=2$, $r(2)=4$, $r(3)=3$, and $r(4)=1$). 
Then, the following condition must be satisfied:
\begin{equation}\label{condicion1}
t_{r(i)}+x_{r(i)}(k-1)\leq x_i(k) \;.
\end{equation} 
A second constraint follows from the demand that trains must connect. 
This gives rise to the following condition 
\begin{equation}\label{condicion2}
t_j+x_j(k-1)\leq x_i(k) \;,\;\forall j\in C(i)\; ,
\end{equation} 
where $C(i)$ is the set of indexes of all the directions of the trains which 
have to provide a connection with the train which leaves in direction $d_i$
(in the case of the network given in Figure~\ref{figure1} we have: $C(1)=\emptyset$, 
$C(2)=\left\{3\right\}$, $C(3)=\left\{1,4\right\}$, and $C(4)=\left\{3\right\}$).  
Finally, the last condition is that a train cannot leave before its 
scheduled departure time. This yields 
\begin{equation}\label{condicion3}
u_i(k)\leq x_i(k) \;,
\end{equation}
where $u_i(k)$ denotes the scheduled departure time for the $k$-th 
train in direction $d_i$. Now, if we assume that a train leaves as soon as 
all the previous conditions have been satisfied, in max-plus notation 
conditions~\eqref{condicion1}, \eqref{condicion2} and~\eqref{condicion3} 
lead to
\begin{equation}\label{ecuacion1}
x_i(k)=\bigoplus_{j\in C(i)} t_j x_j(k-1)\oplus t_{r(i)} x_{r(i)}(k-1)\oplus u_i(k) \;.
\end{equation} 
Therefore, if we define the matrix $A=(a_{ij})\in \zmax^{n\times n}$ by: 
\[
a_{ij}= \left\{
\begin{array}{ll}
t_j & \mbox{ if } j \in C(i) \cup \{r(i)\} , \\
-\infty & \mbox{ otherwise,}
\end{array}
\right.
\]    
then~\eqref{ecuacion1} can be written in matrix form as
\begin{equation}\label{ecuacion2}
x(k)=A x(k-1) \oplus u(k) \;, 
\end{equation} 
where $x(k)=(x_1(k),\ldots ,x_n(k))^T$ and $u(k)=(u_1(k),\ldots ,u_n(k))^T$, which is a system of 
the form of~\eqref{dynamicsystem}. In the particular case of the railway network shown in 
Figure~\ref{figure1} we have 
\[
A=
\begin{pmatrix}
 -\infty & 17 & -\infty & -\infty \\ 
-\infty  & -\infty & 11 & 9 \\
14 & -\infty & 11 & 9 \\
14 & -\infty & 11 & -\infty 
\end{pmatrix}\; .
\] 

Suppose now that we want to decide whether there exists a timetable 
such that the time between two consecutive train departures  
in the same direction is less than a certain given bound or such 
that the time that passengers have to wait to make some connections 
is less than another given bound. To be able to model this kind of 
requirement it is convenient to introduce the extended state vector 
$\overline{x}(k)=(x_1(k),\ldots ,x_n(k),x_1(k-1),\ldots ,x_n(k-1))^T$. 
Then~\eqref{ecuacion2} can be rewritten as 
$\overline{x}(k)=\overline{A}\overline{x}(k-1)\oplus \overline{B}u(k)$, 
where 
\[
\overline{A}=
\begin{pmatrix}
A & \varepsilon \\ 
I  & \varepsilon  
\end{pmatrix}\; \mbox{ and }\;
\overline{B}=
\begin{pmatrix}
I \\ 
\varepsilon  
\end{pmatrix}
\]
(here $I,\varepsilon \in \zmax^{n\times n}$ denote the max-plus 
identity and zero matrices, respectively).   
Assume that we want the time between two consecutive train departures 
in direction $d_i$ to be less than $L_i$ time units. This can be 
expressed as $\overline{x}_i(k)-\overline{x}_{i+n}(k)\leq L_i$, 
or equivalently as $\overline{x}_i(k)-L_i\leq \overline{x}_{i+n}(k)$. 
For simplicity we will take the same bound $L$ for all the directions, 
although everything that follows can be done with different bounds. 
Then the previous condition can be written in matrix form as 
\begin{equation}\label{especificacion1}
\begin{pmatrix}
\varepsilon & \varepsilon \\ 
(-L) I  & \varepsilon  
\end{pmatrix}
\overline{x}(k) \leq \overline{x}(k)\; ,\; \forall k\in \N \; . 
\end{equation}  
Suppose now that we want passengers coming from direction 
$d_i$ not to have to wait more than $M_{ij}$ time units for the 
departure of the train which leaves in direction $d_j$. This can be 
expressed as $\overline{x}_j(k)-a_{ji}-\overline{x}_{i+n}(k)\leq M_{ij}$, 
which is equivalent to 
$\overline{x}_j(k)-a_{ji}-M_{ij}\leq \overline{x}_{i+n}(k)$. Once again, 
if for simplicity we take the same bound $M$ for all the possible 
connections, the previous condition can be written in matrix form as
\begin{equation}\label{especificacion2}
\begin{pmatrix}
\varepsilon & \varepsilon \\ 
(-M) S  & \varepsilon  
\end{pmatrix}
\overline{x}(k) \leq \overline{x}(k)\; , \; \forall k\in \N \; ,
\end{equation}
where the matrix $S=(s_{ij})\in \zmax^{n\times n}$ is defined 
by: $s_{ij}=-a_{ji}$ if $a_{ji}\not = -\infty$ and $s_{ij}= -\infty$ 
otherwise. Finally, in order to have {\em realistic} initial states for the extended state 
vector, we can consider the obvious physical constraints $x(k-1)\leq x(k)$ and 
$Ax(k-1)\leq x(k)$, which lead to the following condition: 
\begin{equation}\label{especificacion3}
\begin{pmatrix}
\varepsilon & I\oplus A \\ 
\varepsilon  & \varepsilon  
\end{pmatrix}
\overline{x}(k) \leq \overline{x}(k)\; , \; \forall k\in \N \; .
\end{equation}  
Therefore, to get the desired behavior of the network, 
the timetable $u(k)$ should be such that the extended state 
vector satisfies conditions~\eqref{especificacion1}, 
\eqref{especificacion2} and~\eqref{especificacion3}, that is, 
such that $E\overline{x}(k) \leq \overline{x}(k)$ for all 
$k\in \N$, 
where 
\[
E=
\begin{pmatrix}
\varepsilon & I\oplus A  \\ 
(-M) S\oplus (-L) I  & \varepsilon  
\end{pmatrix}
\; .
\]
For instance, let us take $L=15$ and $M=4$ in the case of the 
railway network shown in Figure~\ref{figure1}. Then 
$E\overline{x}(k) \leq \overline{x}(k)$ is equivalent to  
$\overline{x}(k)\in \im E^*$ (see the proof of 
Lemma~\ref{lemaHstar}), where 
\[
E^*=
\begin{pmatrix}
 0 & 2 & -2 & -2 & 12 & 17 & 13 & 11 \\ 
-5 & 0 & -4 & -4 & 10 & 12 & 11 & 9 \\
-1 & 1 & 0 & -3 & 14 & 16 & 12 & 10 \\
-1 & 1 & -3 & 0 & 14 & 16 & 12 & 10 \\
-15 & -13 & -17 & -17 & 0 & 2 & -2 & -4 \\
-20 & -15 & -19 & -19 & -5 & 0 & -4 & -6 \\
-16 & -14 & -15 & -15 & -1 & 1 & 0 & -5 \\
-14 & -12 & -13 & -15 & 1 & 3 & -1 & 0 
\end{pmatrix}\; .
\] 
Therefore, our problem is to determine the maximal geometrically 
$(\overline{A},\overline{B})$-invariant semimodule contained 
in $\cK=\im E^*$. With this aim we compute the sequence 
of semimodules $\{\cX_r\}_{r\in \N}$ defined by~\eqref{algoABinv} 
following the method described in Remark~\ref{ObsComputo}  
(which has been implemented with scilab, see~\cite{toolbox}). 
Since the entries of $E^*$ are all finite, from Corollary~\ref{corvolfin} 
we know that this sequence must stabilize. In fact, we have: 
$\cX_5=\cX_4\varsubsetneq \cX_3\varsubsetneq \cX_2\varsubsetneq \cX_1=\cK$. 
Then, the maximal geometrically $(\overline{A},\overline{B})$-invariant 
semimodule $\cK^*$ contained in $\cK$ is $\cX_4$, 
which is generated by the columns of the following matrix
\[
\begin{pmatrix}
17 & 17 & 17 & 18 & 17 \\ 
15 & 15 & 14 & 15 & 15 \\
18 & 18 & 17 & 18 & 18 \\
19 & 19 & 18 & 19 & 19 \\
4  & 2  & 2  & 3  & 2  \\
0  & 0  & 0  & 0  & 0  \\
4  & 4  & 3  & 4  & 4  \\
5  & 5  & 4  & 5  & 2 
\end{pmatrix}\; .
\]
Consequently, it is possible to obtain the desired 
behavior of the network with a suitable choice of 
the timetable $u(k)$ when the initial state belongs 
to $\cK^*$. To be able to compute these 
timetables we use the method described at the end of 
Section~\ref{algABinvSec} to decide whether $\cK^*=\cK_4$ is 
an algebraically $(\overline{A},\overline{B})$-invariant semimodule
(that is, we apply the min-max Howard algorithm to find a state feedback). 
In this way we can see that $\cK^*$ is algebraically 
$(\overline{A},\overline{B})$-invariant and one possible state 
feedback is given by
\[
\overline{F}=
\begin{pmatrix}
14 & 14 & 14 & 13 & 14 & 14 & 14 & 14 \\ 
11 & 14 & 11 & 10 & 14 & 14 & 14 & 14 \\
14 & 14 & 14 & 13 & 14 & 14 & 14 & 14 \\
14 & 14 & 14 & 14 & 14 & 14 & 14 & 14 
\end{pmatrix}\; .
\]
For instance, let us consider the evolution of the railway 
network when the initial state is 
$\overline{x}(0)=(17,15,18,19,4,0,4,5)^T \in \cK^*$ and 
the control $\overline{F}$ is applied. In this case we obtain the 
following trajectory $x(k)$ of the system
\[
\begin{pmatrix}
4 \\ 
0 \\
4 \\
5 
\end{pmatrix},\;
\begin{pmatrix}
17 \\ 
15 \\
18 \\
19 
\end{pmatrix},\;
\begin{pmatrix}
32 \\ 
29 \\
32 \\
33 
\end{pmatrix},\;
\begin{pmatrix}
46 \\ 
43 \\
46\\
47 
\end{pmatrix},\;
\begin{pmatrix}
60 \\ 
57 \\
60 \\
61 
\end{pmatrix},\;
\begin{pmatrix}
74 \\ 
71 \\
74 \\
75 
\end{pmatrix},\;
\ldots
\]
which clearly satisfies the constraints imposed on the network. 
However, if no control is applied, we get the following trajectory starting 
from the same initial state 
\[
\begin{pmatrix}
4 \\ 
0 \\
4 \\
5 
\end{pmatrix},\;
\begin{pmatrix}
17 \\ 
15 \\
18 \\
19 
\end{pmatrix},\;
\begin{pmatrix}
32 \\ 
29 \\
31 \\
31 
\end{pmatrix},\;
\begin{pmatrix}
46 \\ 
42 \\
46\\
46 
\end{pmatrix},\;
\begin{pmatrix}
59 \\ 
57 \\
60 \\
60 
\end{pmatrix},\;
\begin{pmatrix}
74 \\ 
71 \\
73 \\
73 
\end{pmatrix},\;
\ldots
\]
which does not satisfy the constraints imposed on the network, since for 
example the passengers coming from station $S$ on the third train (which 
leaves from station $Q$ in direction $d_4$ at time $31$) will have to wait $6$ 
time units for the next departure of a train in direction $d_3$ toward 
station $R$ (which will take place at time $46$). 

If we want to obtain the desired behavior of the network with a 
periodic timetable, that is with a timetable $u(k)$ of the form 
$u(k)=\lambda^k u$, where $\lambda \in \zmax$ and $u\in \zmax^n$, 
then what we can do is to see if the matrix 
$\overline{A}\oplus \overline{B} \overline{F}$ has an eigenvector in $\cK^*$. 
In this case it can be shown that 
$\overline{x}(0)=(17,14,17,18,3,0,3,4)^T \in \cK^*$ 
is an eigenvector of $\overline{A}\oplus \overline{B} \overline{F}$ 
corresponding to the eigenvalue $\lambda =14$, that is, 
the following equality is satisfied:
\[
(\overline{A}\oplus \overline{B} \overline{F})\overline{x}(0)=
14\overline{x}(0)\; .
\]
Therefore, the periodic timetable
\[ 
u(k)=\overline{F}\overline{x}(k-1)=14^{(k-1)}\overline{F}\overline{x}(0)=
14^{(k+1)}\begin{pmatrix}
3 \\ 
0 \\
3 \\
4 
\end{pmatrix}
\]
leads to the desired behavior of the network when the initial state 
is $\overline{x}(0)$. In other words, one train should leave in each 
direction every $14$ time units but the $k$-th departure time of the 
trains in direction $d_1$ and $d_3$, respectively in direction $d_4$, 
should be scheduled $3$ time units, respectively $4$ time units, 
after the $k$-th scheduled departure time of the train in direction $d_2$. 

Let us finally mention that the computations of the examples presented 
in this paper have been checked using the max-plus toolbox of scilab 
(see~\cite{toolbox}).

\section{Conclusion}
In this paper, the classical concept of $(A,B)$-invariant space is extended to linear 
dynamical systems over the max-plus semiring. This extension presents similar 
difficulties to those encountered in dealing with coefficients in a ring rather than 
coefficients in a field. On the one hand, we show that the classical algorithm for the 
computation of the maximal $(A,B)$-invariant subspace contained in a given space, which is 
generalized to the max-plus algebra framework, need not converge in a finite number of steps. 
However, sufficient conditions for the convergence of this algorithm are given. In particular, 
it is shown that these conditions are satisfied by a class of semimodules of practical interest. 
On the other hand, the existence (which is not guaranteed) and the computation of linear state 
feedbacks are also discussed in the case of finitely generated semimodules. 
Finally, we show that this approach is capable of providing solutions 
to some control problems by considering its application to the study of transportation 
networks which evolve according to a timetable.

\bibliographystyle{alpha} 
 
\bibliography{05-085katz} 

\begin{thebibliography}{CHMSM03}

\bibitem[ALP99]{AssLafPer}
J.~Assan, J.~F. Lafay, and A.~M. Perdon.
\newblock On feedback invariance properties for systems over a principal ideal
  domain.
\newblock {\em IEEE Trans. Automat. Control}, 44(8):1624--1628, 1999.

\bibitem[Ass99]{assan}
J.~Assan.
\newblock {\em Analyse et synth\`ese de l'approche g\'eom\'etrique pour les
  syst\`emes lin\'eaires sur un anneau}.
\newblock Th\`ese de doctorat, Universit\'e de Nantes, France, Octobre 1999.

\bibitem[BCFH99]{BoiCott99}
J.~L. Boimond, B.~Cottenceau, J.~L. Ferrier, and L.~Hardouin.
\newblock Synthesis of greatest linear feedback for timed-event graphs in
  dioid.
\newblock {\em IEEE Trans. Automat. Control}, 44(6):1258--1262, 1999.

\bibitem[BCG03]{bcg99}
P.~Butkovi\v{c} and R.~Cuninghame-Green.
\newblock The equation {$A\otimes x=B\otimes y$} over {$(\mathbb{R}\cup
  \{-\infty\},\max,+)$}.
\newblock {\em Theor. Comp. Sci.}, 293(1):3--12, 2003.

\bibitem[BCOQ92]{bcoq}
F.~Baccelli, G.~Cohen, G.~J. Olsder, and J.~P. Quadrat.
\newblock {\em Synchronization and Linearity}.
\newblock Wiley, Chichester, 1992.

\bibitem[BFHM00]{BoiHar00}
J.~L. Boimond, J.~L. Ferrier, L.~Hardouin, and E.~Menguy.
\newblock Just-in-time control of timed event graphs: update of reference
  input, presence of uncontrollable input.
\newblock {\em IEEE Trans. Automat. Control}, 45(11):2155--2159, 2000.

\bibitem[BH84]{butkovicH}
P.~Butkovi\v{c} and G.~Heged\"{u}s.
\newblock An elimination method for finding all solutions of the system of
  linear equations over an extremal algebra.
\newblock {\em Ekonomicko-matematicky Obzor}, 20(2):203--215, 1984.

\bibitem[BJ72]{BlythJan72}
T.~S. Blyth and M.~F. Janowitz.
\newblock {\em Residuation Theory}.
\newblock Pergamon Press, London, 1972.

\bibitem[BM91]{BasMar91}
G.~Basile and G.~Marro.
\newblock {\em Controlled and Conditioned Invariants in Linear System Theory}.
\newblock Prentice Hall, 1991.

\bibitem[Bra91]{braker91}
J.~G. Braker.
\newblock Max-algebra modelling and analysis of time-dependent transportation
  networks.
\newblock In {\em Proceedings of the 1st European Control Conference}, pages
  1831--1836, Grenoble, France, July 1991.

\bibitem[Bra93]{braker}
J.~G. Braker.
\newblock {\em Algorithms and Applications in Timed Discrete Event Systems}.
\newblock Ph. {D}. {T}hesis, Faculty of Technical Mathematics and Informatics,
  Delft University of Technology, Delft, The Netherlands, 1993.

\bibitem[CDQV85]{cohen85a}
G.~Cohen, D.~Dubois, J.~P. Quadrat, and M.~Viot.
\newblock A linear system theoretic view of discrete event processes and its
  use for performance evaluation in manufacturing.
\newblock {\em IEEE Trans. on Automatic Control}, AC--30:210--220, 1985.

\bibitem[CGQ99]{ccggq99}
G.~Cohen, S.~Gaubert, and J.~P. Quadrat.
\newblock Max-plus algebra and system theory: where we are and where to go now.
\newblock {\em Annual Reviews in Control}, 23:207--219, 1999.

\bibitem[CGQ01]{gaubert01a}
G.~Cohen, S.~Gaubert, and J.~P. Quadrat.
\newblock Duality of idempotent semimodules.
\newblock In {\em Proceedings of the Satellite Workshop on Max-Plus Algebras,
  IFAC SSSC'01}, Praha, 2001. Elsevier.

\bibitem[CH83]{CassaHo83}
C.~G. Cassandras and Y.-C. Ho.
\newblock A new approach to the analysis of discrete event dynamic systems.
\newblock {\em Automatica J. IFAC}, 19(2):149--167, 1983.

\bibitem[CHMSM03]{CottHar03}
B.~Cottenceau, L.~Hardouin, C.~A. Maia, and R.~Santos-Mendes.
\newblock Optimal closed-loop control of timed event graphs in dioids.
\newblock {\em IEEE Trans. Automat. Control}, 48(12):2284--2287, 2003.

\bibitem[CLO95]{CassLafoOlsd}
C.~G. Cassandras, S.~Lafortune, and G.~J. Olsder.
\newblock Introduction to the modelling, control and optimization of discrete
  event systems.
\newblock In {\em Trends in control (Rome, 1995)}, pages 217--291. Springer,
  Berlin, 1995.

\bibitem[CMQV89]{cohen89a}
G.~Cohen, P.~Moller, J.~P. Quadrat, and M.~Viot.
\newblock Algebraic tools for the performance evaluation of discrete event
  systems.
\newblock {\em IEEE Proceedings: Special issue on Discrete Event Systems},
  77(1):39--58, Jan. 1989.

\bibitem[CP94]{conte94}
G.~Conte and A.~M. Perdon.
\newblock Problems and results in a geometric approach to the theory of systems
  over rings.
\newblock In {\em Linear algebra for control theory}, volume~62 of {\em IMA
  Vol. Math. Appl.}, pages 61--74. Springer, New York, 1994.

\bibitem[CP95]{conte95}
G.~Conte and A.~M. Perdon.
\newblock The disturbance decoupling problem for systems over a ring.
\newblock {\em SIAM J. Control Opt.}, 33(3):750--764, 1995.

\bibitem[CTGG99]{gg}
J.~Cochet-Terrasson, S.~Gaubert, and J.~Gunawardena.
\newblock A constructive fixed point theorem for min-max functions.
\newblock {\em Dynamics and Stability of Systems}, 14(4):407--433, 1999.

\bibitem[dDD98]{deVDeSdeM98}
R.~{de Vries}, B.~{De Schutter}, and B.~{De Moor}.
\newblock On max-algebraic models for transportation networks.
\newblock In {\em Proceedings of the International Workshop on Discrete Event
  Systems (WODES'98)}, pages 457--462, Cagliari, Italy, August 1998.

\bibitem[Gau92]{gaubert92a}
S.~Gaubert.
\newblock {\em Th\'eorie des syst\`emes lin\'eaires dans les dio\"\i des}.
\newblock Th\`ese, \'Ecole des Mines de Paris, France, July 1992.

\bibitem[Gau98]{gaubert98n}
S.~Gaubert.
\newblock Exotic semirings: Examples and general results.
\newblock Support de cours de la 26$^{\mbox{\rm i\`eme}}$ \'Ecole de Printemps\
  d'Informatique Th\'eorique, Noirmoutier, 1998.

\bibitem[GG98]{cras}
S.~Gaubert and J.~Gunawardena.
\newblock The duality theorem for min-max functions.
\newblock {\em C.R. Acad. Sci.}, 326(1):43--48, 1998.

\bibitem[GK95]{GargKumar95}
V.~K. Garg and R.~Kumar.
\newblock {\em Modeling and Control of Logical Discrete Event Systems}.
\newblock Kluwer Academic Publisher, Norwell Massachusetts, 1995.

\bibitem[GK03]{gk03}
S.~Gaubert and R.~D. Katz.
\newblock Reachability and invariance problems in max-plus algebra.
\newblock In L.~Benvenuti, A.~{De Santis}, and L.~Farina, editors, {\em
  Proceedings of POSTA'03}, number 294 in Lecture Notes in Control and Inf.
  Sci., pages 15--22, Berlin, Aug. 2003. Springer.

\bibitem[GK04]{gk02a}
S.~Gaubert and R.~D. Katz.
\newblock Rational semimodules over the max-plus semiring and geometric
  approach to discrete event systems.
\newblock {\em Kybernetika}, 40(2):153--180, 2004.
\newblock Also e-print arXiv:math.OC/0208014.

\bibitem[GP97]{maxplus97}
S.~Gaubert and M.~Plus.
\newblock Methods and applications of {$(\max,+)$} linear algebra.
\newblock In R.~Reischuk and M.~Morvan, editors, {\em 14th Symposium on
  Theoretical Aspects of Computer Science, STACS 97 (L\"ubeck)}, volume 1200 of
  {\em Lecture Notes in Comput. Sci.}, pages 261--282, Berlin, 1997. Springer.

\bibitem[Hau82]{hautus82}
M.~L.~J. Hautus.
\newblock Controlled invariance in systems over rings.
\newblock In {\em Feedback control of linear and nonlinear systems
  (Bielefeld/Rome, 1981)}, volume~39 of {\em Lecture Notes in Control and
  Inform. Sci.}, pages 107--122. Springer, Berlin, 1982.

\bibitem[Hau84]{hautus84}
M.~L.~J. Hautus.
\newblock Disturbance rejection for systems over rings.
\newblock In {\em Mathematical theory of networks and systems (Beer Sheva,
  1983)}, volume~58 of {\em Lecture Notes in Control and Inform. Sci.}, pages
  421--432. Springer, London, 1984.

\bibitem[Kat03]{katz}
R.~D. Katz.
\newblock {\em Problemas de alcanzabilidad e invariancia en el \'algebra
  max-plus}.
\newblock Ph. {D}. {T}hesis, Universidad Nacional de Rosario, Argentina,
  November 2003.

\bibitem[Kli03]{klimann99}
I.~Klimann.
\newblock A solution to the problem of {$(A,B)$}-invariance for series.
\newblock {\em Theoret. Comput. Sci.}, 293(1):115--139, 2003.

\bibitem[Lho03]{Lhommeau}
M.~Lhommeau.
\newblock {\em \'{E}tude de syst\`emes \`a \'ev\'enements discrets dans
  l'alg\`ebre $(\max,+)$}.
\newblock Th\`ese de doctorat, ISTIA - Universit\'e d'Angers, France, December
  2003.

\bibitem[LT01]{leboudec}
J-Y. {Le Boudec} and P.~Thiran.
\newblock {\em Network calculus}.
\newblock Number 2050 in LNCS. Springer, 2001.

\bibitem[OSG98]{OlsSubGett98}
G.~J. Olsder, S.~Subiono, and M.~Mc Gettrick.
\newblock On large scale max-plus algebra model in railway systems.
\newblock In {\em Proceedings of the International Workshop on Discrete Event
  Systems (WODES'98)}, Cagliari, Italy, August 1998.

\bibitem[Pin98]{pin95}
J-E. Pin.
\newblock Tropical semirings.
\newblock In J.~Gunawardena, editor, {\em Idempotency (Bristol, 1994)},
  volume~11 of {\em Publ. Newton Inst.}, pages 50--69. Cambridge University
  Press, Cambridge, 1998.

\bibitem[Plu98]{toolbox}
M.~Plus.
\newblock Max-plus toolbox of scilab.
\newblock Available from the contrib section of
  http://www-rocq.inria.fr/scilab, 1998.

\bibitem[RW87]{ramadge87a}
P.~J. Ramadge and W.~M. Wonham.
\newblock Supervisory control of a class of discrete event processes.
\newblock {\em SIAM J. Control and Optimization}, 25(1):206--230, Jan 1987.

\bibitem[WB98]{walkup}
E.~A. Walkup and G.~Borriello.
\newblock A general linear max-plus solution technique.
\newblock In J.~Gunawardena, editor, {\em Idempotency (Bristol, 1994)},
  volume~11 of {\em Publ. Newton Inst.}, pages 406--415. Cambridge University
  Press, Cambridge, 1998.

\bibitem[Won85]{wonham}
W.~M. Wonham.
\newblock {\em Linear multivariable control: a geometric approach}.
\newblock Springer, 1985.
\newblock Third edition.

\end{thebibliography}

\end{document}